\documentclass[11pt]{article}
\usepackage{graphicx}
\usepackage{pdfsync}
\usepackage{psfrag}
\usepackage{amsmath}
\usepackage{amssymb}
\usepackage{epstopdf}
\usepackage{mathdots}
\usepackage{relsize}
\usepackage{color}

\newtheorem{theorem}{Theorem}
\newtheorem{proposition}[theorem]{Proposition}
\newtheorem{lemma}[theorem]{Lemma}
\newtheorem{corollary}[theorem]{Corollary}
\newtheorem{definition}[theorem]{Definition}
\newtheorem{conjecture}[theorem]{Conjecture}

\title{Twisted Alexander polynomials of  2-bridge knots}
\author{
Jim Hoste\\
Pitzer College\\
\\
Patrick D. Shanahan\\
Loyola Marymount University
}
\begin{document}
\maketitle

\begin{abstract} We investigate the twisted Alexander polynomial of a 2-bridge knot associated to a Fox coloring. For several  families of 2-bridge knots, including but not limited to,  torus knots and genus-one knots, we derive  formulae for these twisted Alexander polynomials.  We  use these formulae to confirm a conjecture of Hirasawa and Murasugi for these knots.

\end{abstract}
\section{Introduction} Associated to every knot $K$ is its classical Alexander polynomial $\Delta_K(t)$. This well-known polynomial invariant encodes metabelian information of the knot group. That is, one can think of the Alexander polynomial as describing the representations of the knot group into a metabelian subgroup of a linear group.  Knot groups, however, allow for many other types of non-abelian linear representations. From this perspective, the Alexander polynomial is just a special case of a more general invariant called the twisted Alexander polynomial introduced by Xiao-Song Lin \cite{Lin:2001}. The twisted Alexander polynomial  $\tilde \Delta_{K, \phi}(t)$ depends not only on the knot $K$ but also on a choice of some linear representation $\phi: \pi_1(\mathbb S^3-K) \to GL(n, \mathbb C)$ and it has been studied extensively in recent years.  

In this article, we study  twisted Alexander polynomials of 2-bridge knots associated with a linear representation coming from a Fox coloring of the knot.  If $K_{p/q}$ is a 2-bridge knot and $\ell$ is an odd prime dividing the determinant $q$ of $K_{p/q}$, then it is well-known that there is a Fox $\ell$-coloring of the knot. This coloring provides a linear representation of the knot group that factors through the dihedral group.  We  call the twisted Alexander polynomial associated with this linear representation the {\it $\ell$-twisted} Alexander polynomial of $K_{p/q}$ and denote it by $\tilde \Delta_{p/q}^\ell(t)$. Our study of these polynomials was motivated by the following conjecture of Hirasawa and Murasugi \cite{HM:2009}.

\begin{conjecture}\label{conjecture}\mbox{\rm [Hirasawa and Murasugi]}\label{HM conjecture} Let $K$ be a knot and $\phi: \pi_1(S^3-K)\to D_\ell \to GL(\ell, \mathbb C)$ a non-abelian linear representation, where $D_\ell$ is the dihedral group of order $2\ell$ with $\ell$  an odd prime. Then
$$ \tilde \Delta_{K, \phi}(t)=\frac{\Delta_K(t)}{t-1}f(t)f(-t)$$ for some integer Laurent polynomial $f(t)$ and furthermore
$$f(t)\equiv \left (\frac{\Delta_K(t)}{t+1}\right )^\frac{\ell-1}{2} \ (\mbox{\rm mod } \ell).$$
\end{conjecture}

Hirasawa and Murasugi verify this conjecture for 2-bridge knots whose fundamental groups map onto certain free products but do so without finding a specific formula for the twisted Alexander polynomial. We will extend their work but employ a different method in this paper.  In particular, we will derive a specific formula for the twisted Alexander polynomial of the knots we study and then use it to verify the conjecture directly.  Before we proceed, however, it is important to note that there is some potential ambiguity in the statement of the conjecture.  First of all, neither the Alexander nor the twisted Alexander polynomial is well-defined. The former is only well-defined up to multiples of $\pm t$ while the latter up to multiples of $\pm t^{\ell}$. Thus, the Laurent polynomial $f(t)$ in the conjecture also fails to be well-defined. Because of this there are two different interpretations of the conjecture. A weak form of the conjecture would allow for some  choices of $ \tilde \Delta_{K, \phi}(t)$,  $\Delta_K(t)$, and  $f(t)$ that give equality in the first part and possibly different choices of $\Delta_K(t)$ and  $f(t)$ that give congruency  in the second part.  A stronger form of the conjecture would require the same choices in both parts. For the families of knots studied in this paper, we will prove the stronger version of the conjecture.

Motivated by a beautiful formula for $\Delta_{K_{p/q}}(t)$, we develop a computational method and accompanying formula for $\tilde \Delta_{p/q}^\ell(t)$. Combining this formula with a recursive technique, we  derive specific formulae for several bi-infinite families of 2-bridge knots. In particular, given a specific choice of $p$, $q$, and $\ell$, we derive a formula for the $\ell$-twisted Alexander polynomial of the 2-bridge knot  associated to the fraction 
$$\frac{p'}{q'}=\frac{p+2 \ell j r}{q+2\ell(kp+jar+2 \ell j k r)},$$ 
where $k$ and $j$ are any non-negative integers and $q=ap+r$ with $0<r<a$.  For any 2-bridge knot, the Alexander polynomial is always a factor of its twisted Alexander polynomial and, for these bi-infinite families, the formula for 
$\tilde \Delta_{p'/q'}^\ell(t)/\Delta_{K_{p'/q'}}(t)$
 depends on $p,\ q,\ \ell$, and $k$ but   not on  $j$. 
With this formula in hand, we can verify Conjecture~\ref{conjecture} for every knot in the family. From a practical point of view, one can compute the $\ell$-twisted Alexander polynomial for as many bi-infinite families of this kind as one likes. Furthermore, if we fix $p$ and $\ell$, then only a finite number of computations are needed to derive a formula for  $\tilde \Delta_{p/q}^\ell(t)$ for all possible $q$. We do this for a small sample of $(p, \ell)$, namely: $(5,3)$, $(7,3)$, $(3, 5)$, and $(3, 7)$.

The paper will proceed as follows. In Section~2, we define the $\ell$-twisted Alexander polynomial and find a formula for it which is derived from a simple graph associated to $p/q$ which we call the epsilon graph.  In Section~3, we establish several properties of this graph which are used in Section~4 to determine how $\tilde \Delta_{p/q}^\ell(t)$ changes when $p$ and $q$ are changed in certain ways. Finally, in Section~5, we apply these results to establish Conjecture~\ref{HM conjecture} for several infinite families of 2-bridge knots.

\section{The $\ell$-twisted Alexander polynomial}\label{twisted alex poly}

Consider the 2-bridge knot $K_{p/q}$ determined by the integers $p$ and $q$ where $0<p<q$, both $p$ and $q$ are odd, and $\gcd(p,q)=1$. It is well known that the fundamental group is presented as
$$G_{p/q}=<a, b\, |\, aw=wb>$$
where $w=b^{\epsilon_1}a^{\epsilon_2}\dots a^{\epsilon_{q-1}}$ and $\epsilon_i=(-1)^{\lfloor ip/q \rfloor}$. Here $\lfloor x \rfloor$ denotes the greatest integer less than or equal to $x$. As an aside,  it is easy to show that the $\epsilon_i$'s are symmetric, that is, $\epsilon_i=\epsilon_{q-i}$. 

The Alexander polynomial can be derived from this presentation using the Fox free differential calculus \cite{CF}. Letting $\cal R$ be the relator ${\cal R}=awb^{-1}w^{-1}$ in $G_{p/q}$, a calculation shows that the matrix ${\cal F}=\left( \frac{\partial {\cal R}}{\partial a} \ \  \frac{\partial {\cal R}}{\partial b} \right )$ of Fox derivatives is 
$${\cal F}=\left( 1+(a-1)\frac{\partial w}{\partial a} \ \ \ \  (a-1)\frac{\partial w}{\partial b}-w\right)$$ 
where
\begin{equation}\label{partial of r wrt a} 
\frac{\partial w}{\partial a}=\epsilon_2 b^{\epsilon_1}a^{\frac{\epsilon_2-1}{2}}
+ \epsilon_4 b^{\epsilon_1} a^{\epsilon_2} b^{\epsilon_3}a^{\frac{\epsilon_4-1}{2}}+\dots +
 \epsilon_{q-1} b^{\epsilon_1} a^{\epsilon_2} b^{\epsilon_3}\dots b^{\epsilon_{q-2}}a^{\frac{\epsilon_{q-1}-1}{2}}
\end{equation}
and
\begin{equation}\label{partial of r wrt b}
\frac{\partial w}{\partial b}=\epsilon_1 b^{\frac{\epsilon_1-1}{2}}
+ \epsilon_3 b^{\epsilon_1} a^{\epsilon_2} b^{\frac{\epsilon_3-1}{2}}+\dots +
 \epsilon_{q-2} b^{\epsilon_1} a^{\epsilon_2} b^{\epsilon_3}\dots b^{\frac{\epsilon_{q-2}-1}{2}}.
\end{equation}

To find the Alexander polynomial, we replace $a$ and $b$ in $\cal F$ with their images in the abelianization of $G_{p/q}$, whose generator we call $t$, to obtain the $1\times2$ matrix $A$. In this case both $a$ and $b$ have image $t$. Now  striking either entry of $A$ and taking the determinant of the remaining  $1 \times 1$ matrix gives the following formula for the Alexander polynomial of an arbitrary 2-bridge knot. (See \cite{Minkus:1982} and compare \cite{Hartley:1979}).
\begin{equation}\label{alex poly formula}
\Delta_{K_{p/q}}(t)=1-t^{\epsilon_1}+t^{\epsilon_1+\epsilon_2}-t^{\epsilon_1+\epsilon_2+\epsilon_3}+\dots+t^{\epsilon_1+\epsilon_2+\epsilon_3+\dots +\epsilon_{q-1}}.
\end{equation}
For any 2-bridge knot, we can use (\ref{alex poly formula}) to produce a canonical choice for the Alexander polynomial and we shall do so unless otherwise noted throughout the rest of the paper. Formula~(\ref{alex poly formula}) can also be interpreted graphically in a beautiful way using what we call the {\it epsilon} graph, $E(p,q)$. (A similar graph is studied by Hirasawa and Murasugi in \cite{HM:2007}.)
\begin{definition}\rm
Given $0<p<q$, $pq$ odd, and $\gcd(p,q)=1$,  the {\it epsilon graph}, $E(p, q)$, is the set of  $q$ vertices  $\{v_i\}_{i=0}^{q-1}$ in the plane whose coordinates are $v_i=(i, \sum_{j=1}^i \epsilon_j)$, together with the set of $q-1$ edges that connect consecutive vertices. The $i$-th edge has slope $\epsilon_i$. A maximal set of consecutive edges in $E(p, q)$ having the same slope is called a {\it segment}.
\end{definition}

\begin{figure}[h]
    \begin{center}
    \leavevmode
    \scalebox{.65}{\includegraphics{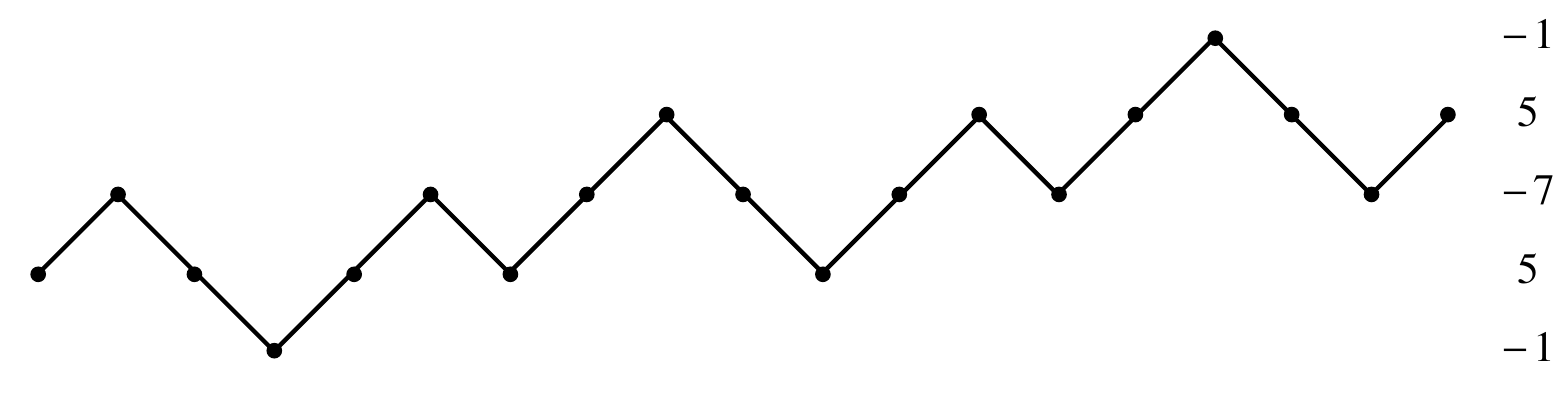}}
    \end{center}
\caption{The epsilon graph of the 2-bridge knot $K_{11/19}$.}
\label{sawgraph}
\end{figure}

We illustrate the case of $p/q=11/19$ in Figure~\ref{sawgraph}.  Summing the number of vertices at each horizontal level of $E(11, 19)$ gives the alternating coefficients of the Alexander polynomial.  That is, $\Delta_{K_{11/19}}(t) = -t^{-1}+5-7 t+5 t^2-t^3$. 
This method of deriving  $\Delta_{K_{p/q}}(t)$ from $E(p,q)$ follows directly from (\ref{alex poly formula}).

 The twisted Alexander polynomial  can also be obtained using the Fox differential calculus. See  \cite{Lin:2001} and \cite{W:1994}. We first choose a linear representation
$\phi: G_{p/q} \to GL(n, \mathbb C)$. The {\it twisted Alexander matrix} $\hat A_{\phi}$ is obtained by replacing $a$, $b$, and $1$ in the first column of $\cal F$ with $t \phi(a)$, $t \phi(b)$,  and the identity matrix $I$, respectively. That is,
\begin{equation}\label{def of ahat}
\hat A_{\phi} = I+(t \phi(a)-I)  \left. \frac{\partial w}{\partial a} \right|_{a= t \phi(a), b=t \phi(b)}.
 \end{equation}

The {\it twisted Alexander polynomial with respect to $\phi$} is then defined (see \cite{W:1994}) as
\begin{equation}
\label{twisted alex formula}
\tilde \Delta_{K, \phi}(t)=\frac{\det {\hat A_{\phi}}}{\det(t\phi(b)-I)}.
\end{equation}

Furthermore, if  $\phi$ is unimodular, this is an invariant of $K$ up to multiples of $\pm t^n$.

We now focus on linear representations obtained from Fox colorings. If $\ell$ is an odd prime, then a {\it Fox coloring}, or an {\it $\ell$-coloring}, is a homomorphsim $\rho$ from $G_{p/q}$ onto the dihedral group $D_\ell$ of $2\ell$ elements.  Recall that the dihedral group has a presentation 
$$D_\ell=<x, y\,| \, x^2=1, y^\ell=1, yx=xy^{-1}>.$$
(Note that the presentation  used in \cite{HM:2009}  uses generators $y^{-1}x y$ and $xy$ in place of $x$ and $y$, respectively.) 
It is well known that  the determinant of the 2-bridge knot $K_{p/q}$ is $q$ and that an $\ell$-coloring exists if and only if $\ell$ divides $q$. Furthermore, if an $\ell$-coloring exists, then, up to conjugation, we may assume that $a$ and $b$ are sent to $x$ and $xy^m$, respectively, for some $0<m<\ell$. Since the twisted Alexander polynomial is defined using determinants, this conjugation will not change $\tilde \Delta_{K, \phi}(t)$.
We now choose, and fix, the following faithful, linear representation $\pi:D_{\ell} \rightarrow GL(\ell, \mathbb C)$ defined by 
\begin{equation}
\label{defn of F R}
\pi(x) = F =\begin{pmatrix} 0 & \dots & 0 & 1 \\ 0 & \dots & 1 & 0 \\ \vdots & \iddots & \vdots & \vdots \\ 1 & \dots & 0 & 0\end{pmatrix} \ \ \ \mbox{and} \ \ \ \pi(y) = R =\begin{pmatrix} 0 & \dots & 0 & 1 \\ 1 & \dots & 0 & 0 \\ \vdots & \ddots & \vdots & \vdots \\ 0 & \dots & 1 & 0\end{pmatrix}.
\end{equation}
Observe that $F$ and $R$ are obtained from the identity matrix by reversing the rows and cyclically permuting the rows, respectively. Thus, an $\ell$-coloring $\rho$ determines a linear representation 
$$\phi = \pi \circ \rho : G_{p/q} \rightarrow GL(\ell, \mathbb C),$$
and hence  an associated twisted Alexander polynomial.

To now determine the twisted Alexander matrix $\hat A_{\phi}$, we replace $\phi(a)$ and $\phi(b)$ in (\ref{def of ahat}) with the matrices $F$ and $FR^m$, respectively. Using the formula for $\partial w/\partial a$ in (\ref{partial of r wrt a}) and the fact that both $F$ and $FR^m$ are their own inverses, we see that 
\begin{equation}\label{alex matrix}{\hat A_{\phi}}=I+\sum_{k=1}^{(q-1)/2}t^{\sum_{j=1}^{2k}\epsilon_j}(I-t^{-\epsilon_{2k}}F)R^{mk \epsilon_{2k}}.
\end{equation}

The numerator of the twisted Alexander polynomial is the determinant of $\hat A_{\phi}$. By examining the form of $FR^m$, it follows that the denominator is
\begin{equation}\label{denominator formula}
\det(tFR^m-I)=-(1+t)^\frac{\ell-1}{2}(1-t)^\frac{\ell+1}{2}.
\end{equation}

Our strategy for finding the twisted Alexander polynomial is to compute the determinant of a similar matrix $P^T \hat A_{\phi} P$ which has a particularly nice form. If  $mn\equiv 1$ (mod $\ell$) and $\lambda=e^{2 \pi i/\ell}$, then we define $P$ to be the  $\ell \times \ell$ symmetric matrix given by $P_{ij}=\lambda^{-(i-1)(j-1)n}$. Now from (\ref{alex matrix}), we see that in order to compute the determinant of $P^T \hat A_{\phi} P=P \hat A_{\phi} P$ we need to determine the matrices $P R^{mk} P$ and $P F R^{mk} P$.

If  $\vec v_j=(1, \lambda^{-jn}, \lambda^{-2 jn}, \dots, \lambda^{-(\ell-1)jn})^T$ for $0 \leq j < \ell$, then a  simple calculation shows that
$$R^m \vec v_j=\lambda^j \vec v_j.$$
Thus, $\lambda^j $ is an eigenvalue of $R^m$ with associated eigenvector $\vec v_j$. Furthermore, note that the $j$-th column of $P$ is the eigenvector $\vec v_{j-1}$.  Using the definition of $\lambda$, we see that 
$$\vec v_i\cdot \vec v_j=\begin{cases}
 \ell& \text{for $i+j\equiv 0$ (mod $\ell$)} \\
  0 & \text{otherwise.}  \\
 \end{cases}$$
Letting $\vec e_i$ denote the $i$-th standard basis vector,  we then obtain 
\begin{eqnarray*}
(PR^{mk}P)_{ij}&=&\vec e_i \cdot (P R^{mk}P \vec e_j)\\
&=&\vec e_i \cdot (PR^{mk} \vec v_{j-1})\\
&=&\lambda^{k(j-1)} \, (\vec e_i^{\ T} P) \vec v_{j-1}\\
&=&\lambda^{k(j-1)} \, \vec v_{i-1} \cdot \vec v_{j-1}\\
&=&\begin{cases}
\ell \lambda^{k(j-1)} &\text{if $i+j \equiv 2$ (mod $\ell$),}\\
0&\text{otherwise.}
\end{cases}
\end{eqnarray*}
Therefore,
$$P R^{mk}P=\ell
\begin{pmatrix}
1&0&0&\cdots &0&0\\
0&0&0&\cdots&0&\lambda^{(\ell-1)k}\\
0&0&0&\cdots&\lambda^{(\ell-2)k}&0\\
\vdots&\vdots&\vdots&\ddots&\vdots&\vdots\\
0&0&\lambda^{2k}&\cdots&0&0\\
0&\lambda^k&0&\cdots&0&0
\end{pmatrix}.$$

Furthermore, if we set $k=0$, then we find that $\det (P P)=(-1)^\frac{\ell-1}{2}\ell^{\ell}$.

Consider now the product $P F R^{mk}P$. Then
\begin{align*}
(P F R^{mk}P)_{ij} &= \vec e_i \cdot (P F R^{mk} P \vec e_j) \\
& = ( \vec e_i^{\ T} P) F R^{mk} \vec v_{j-1}\\
&= \lambda^{k(j-1)} \, \vec v_{i-1}^{\ T} \, F \, \vec v_{j-1}\\
&=\lambda^{k(j-1) }(\lambda^{-(\ell-1)(i-1)n}, \ldots, \lambda^{-2 (i-1)n}, \lambda^{-(i-1)n}, 1) \cdot \vec v_{j-1}\\
&=\lambda^{(k+n)(j-1)} \sum_{l=1}^{\ell} \left(\lambda^{(i-j)n}\right)^{l}\\
&=\begin{cases}
\ell \lambda^{(k+n)(j-1)}&\text{if $i=j$,}\\
0&\text{if $i\ne j$.}
\end{cases}
\end{align*}

Therefore,
$$P  F R^{mk}P=\ell
\begin{pmatrix}
1&0&0&\cdots &0&0\\
0&\lambda^{k+n}&0&\cdots&0&0\\
0&0&\lambda^{2(k+n)}&\cdots&0&0\\
\vdots&\vdots&\vdots&\ddots&\vdots&\vdots\\
0&0&0&\cdots&\lambda^{(\ell-2)(k+n)}&0\\
0&0&0&\cdots&0&\lambda^{(\ell-1)(k+n)}
\end{pmatrix}.$$

Combining all of the above, we see that  $P\hat A_{\phi} P$ has the following form:
$$P\hat A_{\phi} P=\ell\begin{pmatrix}
\Delta_{K_{p/q}}&0&0&\cdots&0&0\\
0&d_1&0&\cdots&0&e_{\ell-1}\\
0&0&d_2&\cdots&e_{\ell-2}&0\\
\vdots&\vdots&\vdots&\ddots&\vdots&\vdots\\
0&0&e_2&\cdots&d_{\ell-2}&0\\
0&e_1&0&\cdots&0&d_{\ell-1}
\end{pmatrix}
$$
where $d_i$ and $e_i$ are each Laurent polynomials in $t$ with coefficients in $\mathbb Z[\lambda]$ given by
\begin{equation}\label{formula for d}d_i(\lambda, t)=- \lambda^{in} \sum_{k=1}^{\frac{q-1}{2}}t^{\sum_{j=1}^{2k-1}\epsilon_j}  \lambda^{ik\epsilon_{2k}} \ \ \ \ \ \mbox{and} \ \ \ \ \ e_i(\lambda, t)=1+\sum_{k=1}^{\frac{q-1}{2}}t^{\Sigma_{j=1}^{2k}\epsilon_j}\lambda^{ik\epsilon_{2k}}.\end{equation}
The determinant of $\hat A_\phi$ is now given by
\begin{eqnarray*}
\mbox{det\,}\hat A_{\phi}&=&\frac{\mbox{det}\left(P \hat A_{\phi} P\right)} {\mbox{det}\left(P P\right)}\\
&=& (-1)^\frac{\ell-1}{2} \Delta_{K_{p/q}}(t) \prod_{i=1}^\frac{\ell-1}{2} \Big (d_i(\lambda, t) d_{\ell-i}(\lambda, t) - e_i(\lambda, t) e_{\ell-i}(\lambda, t)\Big).
\end{eqnarray*}

It is important to note that $d_i(\lambda, t)=d_1(\lambda^i, t)$ and $e_i(\lambda, t)=e_1(\lambda^i, t)$. Hence  we will now write  $d_{p/q}(\lambda, t)$ and $e_{p/q}(\lambda, t)$ in place of $d_1(\lambda, t)$ and $e_1(\lambda, t)$, respectively, to emphasize the dependance on $p$ and $q$.  
Moreover, since the factors $-\lambda^{in}$ and $-\lambda^{-in}$ cancel in the product $d_{p/q}(\lambda^i,t)d_{p/q}(\lambda^{-i},t)$ for each $i$,  the determinant of $ \hat A_{\phi} $ and therefore the twisted Alexander polynomial, $\tilde \Delta_{K_{p/q}, \phi}(t)$,  does not depend on $n$ and $m$. In other words, there is no dependence on the choice of homomorphism of the knot group onto the dihedral group. Thus, we will henceforth assume that $m=n=1$. The following definition fixes notation and terminology that we use throughout the rest of the paper.

\begin{definition}\rm Let $K_{p/q}$ be a 2-bridge knot, $\ell$ an odd prime dividing $q$, $\lambda = e^{2 \pi i/\ell}$, $\rho : G_{p/q} \rightarrow D_{\ell}$ any surjective representation, and $\phi = \pi \circ \rho$ where $\pi$ is the faithful linear representation defined in (\ref{defn of F R}).  The  twisted Alexander polynomial of $K_{p/q}$ with respect to $\phi$  will be called the {\it $\ell$-twisted Alexander polynomial} of $K_{p/q}$ and denoted by
$\tilde \Delta_{p/q}^{\ell}(t)$.
\end{definition}

The preceding discussion proves the following theorem.  In this theorem, we use $\doteq$ to denote equality of Laurent polynomials up to factors of the form $\pm t^{r\ell}$ where $r \in \mathbb Z$.

\begin{theorem} \label{twisted alex thm}
If $K_{p/q}$ is a 2-bridge knot and $\ell$ an odd prime dividing $q$, then the $\ell$-twisted Alexander polynomial of $K_{p/q}$ is
\begin{equation}\label{twisted poly formula}
\tilde \Delta_{p/q}^\ell(t) \doteq \frac{\Delta_{K_{p/q}}(t)}{t-1}  \mathlarger{\mathlarger{\prod}}_{i=1}^{\frac{\ell-1}{2}}\frac{ d_{p/q}(\lambda^i, t)d_{p/q}(\lambda^{-i},t)-e_{p/q}(\lambda^i, t)e_{p/q}(\lambda^{-i},t)}{t^2-1}
\end{equation}

where $\Delta_{K_{p/q}}(t)$ is given by (\ref{alex poly formula}), 
\begin{equation} \label{defn of d e}
d_{p/q}(\lambda, t)= \mathlarger{\sum}_{k=1}^{\frac{q-1}{2}}t^{\sum_{j=1}^{2k-1}\epsilon_j}  \lambda^{k\epsilon_{2k}}, \ \ \ \ \ 
\mbox{and} \ \ \ \ \ 
e_{p/q}(\lambda, t)=1+\mathlarger{\sum}_{k=1}^{\frac{q-1}{2}}t^{\sum_{j=1}^{2k}\epsilon_j}\lambda^{k\epsilon_{2k}}.
\end{equation}
\end{theorem}

Since the determinant $d_{p/q}(\lambda, t)d_{p/q}(\lambda^{-1}, t)-e_{p/q}(\lambda, t)e_{p/q}(\lambda^{-1}, t)$  in (\ref{twisted poly formula}) appears frequently in the remainder of the paper, we will denote it as $D_{p/q}(\lambda, t)$. Note that the right hand side of (\ref{twisted poly formula}) gives a canonical form of $\tilde \Delta_{p/q}^\ell(t)$ that we will use throughout the rest of the paper.

For any 2-bridge knot, the formulae for $d$ and $e$ given in (\ref{defn of d e}) can be read  from the epsilon graph after it has been been labeled as follows. Label the vertex $v_0$ with $0$. Next label both vertices $v_{2i-1}$ and $v_{2i}$ with $\pm i$, where the sign is the slope of the segment between these two vertices (which is $\epsilon_{2i}$). Each vertex now contributes a term of the form $t^s \lambda^r$ to either $d$ or $e$ as follows. First, all vertices on the same horizontal level correspond to the same power $t^s$ just as in the case of using the epsilon graph to compute the Alexander polynomial. Odd vertices contribute odd powers of $t$ to $d$, while even vertices contribute even powers of $t$ to $e$. Finally, the value of $r$ at each vertex is the label of that vertex. 
   
\begin{figure}[h]
    \begin{center}
    \leavevmode
    \scalebox{1.2}{\includegraphics{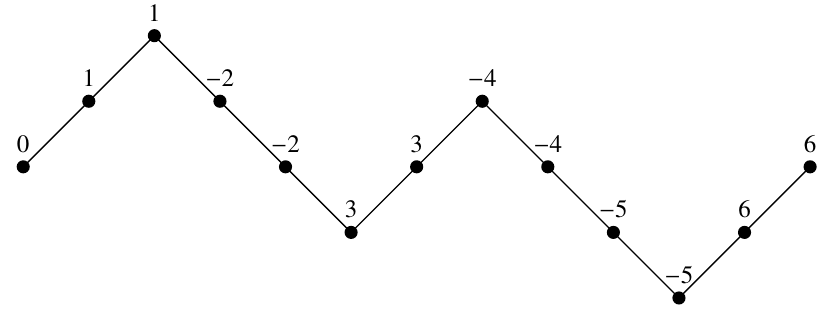}}
    \end{center}
\caption{Labeling the epsilon graph of $K_{5/13}$ so as to read off $d$ and $e$.}
\label{sawtooth labeling}
\end{figure}

We give an example with $p/q=5/13$ in Figure~\ref{sawtooth labeling}. Using the labeling, we obtain
$$e_{5/13}(\lambda,t)=t^{-2}\lambda^{-5}+1+\lambda^{-2}+\lambda^3+\lambda^{-4}+\lambda^6+t^2 \lambda$$ and 
$$d_{5/13}(\lambda,t)=t^{-1}(\lambda^3+\lambda^{-5}+\lambda^6)+t(\lambda+\lambda^{-2}+\lambda^{-4}).$$
Hence
$$
\mathlarger{\mathlarger{\prod}}_{i=1}^6 \frac{D_{5/13}(\lambda^i, t)}{t^2-1}=\frac{(t-1)^6 (t+1)^6
   \left(t^6-2 t^5+t^3-2
   t+1\right)^2 \left(t^6+2
   t^5-t^3+2 t+1\right)^2}{t^{24}}
$$

Now using Theorem~\ref{twisted alex thm}
$$\tilde \Delta_{5/13}^{13}(t)\doteq \frac{\Delta_{K_{5/13}}(t)}{t-1}f(t)f(-t),$$
where $\Delta_{K_{5/13}}(t)=t^{-2}-3t^{-1}+5-3 t+t^2$ and $f(t)=t^{-12}(t+1)^6 (t^6-2 t^5+t^3-2 t+1)^2.$
Notice that the form of $\tilde \Delta_{5/13}^{13}(t)$ is as predicted by Conjecture~\ref{conjecture} and it is easy to verify that the second part of the conjecture is true as well. 

\section{Properties of the Epsilon Graph}

In this section we compile a number of properties of the epsilon graph that will be useful in computing $d$ and $e$ and hence ultimately in computing $\tilde \Delta_{p/q}^\ell(t)$ via Theorem~\ref{twisted alex thm}. Many of the properties in Propositions~\ref{facts 1 about E} and \ref{facts 2 about E} below have already appeared in print, see for example \cite{HM:2007},  but are included here for completeness.

It is nice to think of the  $\epsilon_i$'s in the following way. On the interval from $0$ to $pq$, mark all the points $p, 2p, 3p, \dots , (q-1)p$. If $ip$  lies between $2kq$ and $(2k+1)q$ for some $k$, then $\epsilon_i=1$, while if $ip$ lies between $(2k-1)q$ and $2k q$, then $\epsilon_i=-1$. With this perspective we see that $E(p, q)$ has $p$ segments. This leads naturally to the following definition.
\begin{definition} \rm Given $0<p<q$, $pq$ odd, and $\gcd(p,q)=1$, define the {\it sigma sequence} $\sigma(p, q)=\{\sigma_i\}_{i=1}^p$ by 
letting $\sigma_i$ be the number of edges in the $i$-th segment of $E(p, q)$. Equivalently, $\sigma_i$ is the number of multiples of $p$ contained in the open interval $((i-1)q, iq)$.
\end{definition}

The case of $p=1$ is particularly simple. The epsilon graph consists of a single segment of length $q-1$ because $\epsilon_i=1$ for all $i$.  Hence we will assume $p>1$ for the rest of this section.  We begin with a useful lemma.

\begin{lemma}\label{floor fact}
Suppose $1<p<q$ and $\gcd(p,q)=1$. If $1\le k<p$, then 
$$\left \lfloor \frac{kq}{p}\right \rfloor=\left \lfloor \frac{kq-1}{p}\right \rfloor.$$
\end{lemma}
{\bf Proof:} We may write $kq=mp+r$ where $m>0$ and $1\le r<p$. Now $\frac{kq}{p}=m+\frac{r}{p}$ so $\lfloor \frac{kq}{p}\rfloor=m$. But $kq-1=mp+r-1$ so $\frac{kq-1}{p}=m+\frac{r-1}{p}$ and $0\le \frac{r-1}{p}<\frac{p-1}{p}<1$. Thus $\left \lfloor \frac{kq-1}{p}\right \rfloor=m$ and $\left \lfloor \frac{kq}{p}\right \rfloor=\left \lfloor \frac{kq-1}{p}\right \rfloor$. \hfill $\square$

\begin{proposition}\label{facts 1 about E} Suppose $1<p<q$, $pq$ odd, and $\gcd(p,q)=1$. Further assume that $q=a_1p+r_1$ with  $0<r_1<p$.   Then
\vspace{-.2 in} 
\begin{enumerate}
\item The sigma sequence $\sigma(p, q)$ is symmetric.
\item The graph $E(p, q)$ has a central segment of even length. Moreover, the number of segments of even length is odd and the number of segments of odd length is even.
\item The total number of edges in $E(p,q)$ is $\sum_{i=1}^p\sigma_i=q-1$.
\item For $1\le k\le p$ we have $\sum_{i=1}^k \sigma_i=\left \lfloor \frac{kq-1}{p}\right \rfloor.$ In particular, $\sigma_1=a_1$.
\item For $1\le i\le p$
 we have $\sigma_i=\left \lfloor \frac{iq-1}{p}\right \rfloor-\left \lfloor \frac{(i-1)q-1}{p}\right \rfloor.$
\item The sum of any $k$ consecutive $\sigma_i$'s is equal to either $\sum_{i=1}^k \sigma_i$ or $1+\sum_{i=1}^k \sigma_i$. In particular,  for any $i$, $\sigma_i \in \{\sigma_1, 1+\sigma_1\}$.
\end{enumerate}
\end{proposition}
\noindent{\bf Proof:} Properties 1--3 are immediate.

To prove property~4, which immediately implies property 5, we  look at the length of the interval  $[0, kq]$ in two ways  to obtain
$$p\sum_{i=1}^k \sigma_i+\delta=kq,$$
where $\delta=kq-p\sum_{i=1}^k \sigma_i$ and $0< \delta \le p$. Note that $\delta=p$ if and only if $k=p$. Now
$$\sum_{i=1}^k \sigma_i+\frac{\delta}{p}=\frac{kq}{p}.$$
If $1\le k<p$, then by Lemma~\ref{floor fact}
$$\sum_{i=1}^k \sigma_i=\left \lfloor \frac{kq}{p}\right \rfloor=\left \lfloor \frac{kq-1}{p}\right \rfloor.$$
If instead, $k=p$, then 
 $$\sum_{i=1}^p \sigma_i=q-1=\left \lfloor \frac{pq-1}{p}\right \rfloor.$$

Now consider property 6. The result is trivial if $k=p$, so assume that $1\le k<p$. Consider $k$ consecutive $\sigma_i$'s, say $\sigma_{j+1}, \sigma_{j+2}, \dots, \sigma_{j+k}.$  These  correspond to multiples of $p$ that lie in the open interval $(jq, (j+k)q)$. Let $\alpha$ be the distance between $jq$ and the first multiple of $p$ in $(jq, (j+k)q)$ and let $\beta$ be the distance between $(k+j)q$ and the last multiple of $p$ in $(jq, (j+k)q)$. We have that $0<\alpha \le p$ and $0<\beta \le p$. Now
$$kq=\alpha + p\left (\sum_{i=j+1}^{j+k}\sigma_i-1\right)+\beta.$$
From this we see that   $0<\alpha+\beta<2 p$. Now
$$\frac{kq}{p}=\frac{\alpha+\beta}{p} + \sum_{i=j+1}^{j+k}\sigma_i-1,$$ which gives
$$\left \lfloor \frac{kq}{p}\right \rfloor=\left \lfloor \frac{\alpha+\beta}{p}\right \rfloor + \sum_{i=j+1}^{j+k}\sigma_i-1.$$ 
Since $0<\frac{\alpha+\beta}{p}<2$ the  floor $\left \lfloor \frac{\alpha+\beta}{p}\right \rfloor$ is $0$ or $1$. Hence  $\sum_{i=j+1}^{j+k}\sigma_i$ is equal to either $\left \lfloor \frac{kq}{p}\right \rfloor$ or $1+\left \lfloor \frac{kq}{p}\right \rfloor$. But $\left \lfloor \frac{kq}{p}\right \rfloor=\left \lfloor \frac{kq-1}{p}\right \rfloor=\sum_{i=1}^{k}\sigma_i$ by property~4 and Lemma~\ref{floor fact}.   \hfill $\square$

Proposition~\ref{facts 1 about E} implies that segments come in two possible lengths,  $\sigma_1$ or $1+\sigma_1$. Call segments of length $\sigma_1$ {\it short} and segments of length $1+\sigma_1$ {\it long}. Call a segment {\it isolated} if it is does not have an adjacent segment of the same length. Note that if no long segments exist, then vacuously, all long segments are isolated. By a {\it cluster} of segments we mean   a maximal set of consecutive segments all of the same length. The next proposition lists more facts about $E(p, q)$. 

\begin{proposition}\label{facts 2 about E} Suppose $1<p<q$, $pq$ odd, and $\gcd(p,q)=1$. Assume the Euclidean algorithm gives
$$\begin{array}{rclcl}
q&=&a_1 p+r_1, && 0<r_1<p\\
p&=&a_2 r_1+r_2,&& 0\le r_2<r_1\\
r_1&=&a_3 r_2+r_3,  && 0\le r_3<r_2\\
&\vdots&\\
r_n&=&a_{n+2}r_{n+1}+1
\end{array}
$$
   Then the epsilon graph  has the following properties:
   \vspace{-.2 in}
\begin{enumerate}
\item All short segments are isolated or all long segments are isolated (or both). 
\item The number of long segments is $r_1-1$ and the number of short segments is $p-r_1+1$.
\item The initial cluster of short segments has length  $a_2$.
\item If all long segments are isolated, then every cluster of short segments has length $a_2$ or $a_2-1$. 
\item The number of clusters of short segments of length $a_2$ is equal to $r_2+1$. If $a_2>1$, then  the number of clusters of short segments of length $a_2-1$ equals $r_1-r_2-1$. 
\item If there exists a nonisolated long segment and more than one cluster of long segments, then the first cluster of long segments has length $a_3$ and every cluster of long segments has length $a_3$ or $a_3+1$.
\item If there exists a nonisolated long segment and more than one cluster of long segments, then the number of clusters of long segments of length $a_3$ is $r_2-r_3+1$ and the number of length $a_3+1$ is $r_3-1$.
\end{enumerate}
\end{proposition}
\noindent{\bf Proof:}
If there are two short segments in a row and two long segments in a row, then this would give two pairs of consecutive $\sigma_i$'s whose sums differ by 2, which is not possible by Proposition~\ref{facts 1 about E}. This proves the first property.

To prove property 2, let $L$ be the number of long segments and $p-L$ be the number of short segments. Then 
$$(p-L)\sigma_1+L(\sigma_1+1)=\sum_{i=1}^p\sigma_i=q-1=a_1p+r_1-1.$$
Hence $L=r_1-1$ since $a_1=\sigma_1$.

Property 3 is proven as follows. If there are no long segments, then $r_1=1$ and so $a_2=p$ which is the length of the initial, and only,  cluster of short segments. Alternatively, suppose the initial cluster of short segments has length $k$ and is followed by a long segment. Now $r_1>1$.  This gives the inequalities
$$\begin{array}{rcccl}
k\sigma_1p&<& kq &<& (k\sigma_1+1)p\\
((k+1)\sigma_1+1)p &<&(k+1)q &<&((k+1)\sigma_1+2)p.
\end{array}
$$
If we now replace $q$ with $\sigma_1p+r_1$ we may derive that $k<\frac{p}{r_1}<k+1$ and hence $k=\left \lfloor \frac{p}{r_1}\right \rfloor=a_2$.

One case of property 4 is easy. Suppose some cluster of short segments has length $k\le a_2-2$. By property~4 it cannot be the initial or final cluster, and hence must have  long segments on either end. Now the sum of these $k+2$ consecutive  $\sigma_i$'s is $(k+2)\sigma_1+2$. But, $\sum_{i=1}^{k+2} \sigma_i=(k+2)\sigma_1.$ This contradicts Proposition~\ref{facts 1 about E}.

To prove  the rest of property 4, we begin with the following useful observation. Let $0<\alpha<p$ and consider the set of points $S_\alpha=\{\alpha, \alpha+p, \alpha+2p, \alpha+3p, \dots \}$. It is easy to show that if $0<\alpha<r_1 $, then exactly $\sigma_1+1$ points of $S_\alpha$ lie in the open interval $(0,q)$. On the other hand, if  $r_1 <\alpha<p$, then exactly $\sigma_1$ points of $S_\alpha$ lie in the open interval $(0,q)$. Moreover, in this case, $\alpha+(\sigma_1+1)p-q=\alpha-r_1$. Hence there are exactly $\sigma_1$ points of $S_\alpha$ contained in the open interval $(q, 2q)$ if and only if $\alpha-r_1>r_1$.

Now suppose that there is a cluster of short segments of length more than $a_2$. Suppose the first segment of this cluster lies in the open interval $(jq, (j+1)q)$. Short segments continue at least as far as the interval $((j+a_2)q, (j+a_2+1)q)$. If $\alpha$ is the distance between $jq$ and the first multiple of $p$ in $(jq, (j+1)q)$, then $p>\alpha>r_1$. Furthermore, the distance between $(j+a_2)q$ and the first multiple of $p$ in $((j+a_2)q, (j+a_2+1)q)$ is $\alpha-a_2r_1$ and thus $\alpha-a_2r_1>r_1$.  A similar analysis of the initial cluster of short segments implies that the distance between $a_2 q$ and the first multiple of $p$ in $(a_2q, (a_2+1)q)$ is $p-a_2r_1$. Because the $a_2+1$'st segment is long, we have that $p-a_2r_1<r_1$. These inequalities imply that $\alpha>r_1+a_2 r_1>p$, a contradiction.

We now turn our attention to property 5. Suppose first that $a_2=1$. Let $x$ be the number of clusters of short segments of length 1. Then $x$ equals the number of short segments, so we have $x=p-(r_1-1)$ by property~2. But $p=r_1+r_2$ so we obtain $x=r_2+1$. Suppose now that $a_2>1$.   Let $x$ and $y$ be the number of clusters of short segments of lengths $a_2$ and $a_2-1$, respectively. Since  long segments are isolated, the total number of clusters of short segments is $x+y=r_1$. Counting edges gives
$$q-1=\sigma_1p+r_1-1=xa_2\sigma_1+y(a_2-1)\sigma_1+(r_1-1)(\sigma_1+1)$$
which simplifies to $x=r_2+1$. Furthermore, we have $y=r_1-x=r_1-r_2-1.$

For the first part of property 6, since long segments exist we have $r_1>1$. Furthermore, since long segments are nonisolated, short segments are isolated and so $a_2=1$.  Suppose the first  cluster of long segments has length $k$. This gives the following inequalities.
 $$\begin{array}{rcccl}
((k+1)\sigma_1+k)p&<& (k+1)q &<& ((k+1)\sigma_1+k+1)p\\
((k+2)\sigma_1+k)p&<& (k+2)q &<& ((k+2)\sigma_1+k+1)p.\\
\end{array}
$$
If we replace $q$ with $\sigma_1 p+r_1$ and $p$ with $r_1+r_2$, we obtain $k<\frac{r_1}{r_2}<k+1$ so that $k=\left \lfloor \frac{r_1}{r_2}\right \rfloor=a_3$.

For the second part of property 6, note that the epsilon graph always begins and ends with a cluster of short segments. Hence any cluster of long segments is preceded and followed by a cluster of short segments. Suppose two clusters of long segments have $x$ and $y$ segments respectively with $x+2\le y$. The cluster of length $x$ has short segments on either side. This gives $x+2$ consecutive $\sigma_i$'s that sum to $(x+2)\sigma_1+x$. But the cluster of length $y$ contains a set of $x+2$ consecutive $\sigma_i$'s   that add to $(x+2)(\sigma_1+1)=(x+2)\sigma_1+x+2$. But now the sums of two sets of $x+2$ consecutive $\sigma_i$'s differ by 2, which is impossible by Proposition~\ref{facts 1 about E}. Thus, all  clusters of long segments have length $a_3-1, a_3$, or $a_3+1$. It remains to show that $a_3-1$ is impossible.
Suppose instead that a cluster of long segments exists with length $a_3-1$. Assume the short segment that precedes this cluster lies in the open interval $(jq, (j+1)q)$ and that the short segment that follows it lies in the open interval $((j+a_3)q, (j+a_3+1)q)$. Let $\alpha$ be the distance between $jq$ and the first multiple of $p$ in $(jq, (j+1)q)$. Since a short segment is contained in this interval, we have $r_1<\alpha<p$. The distance between $(j+a_3)q$ and the first multiple of $p$ in the interval $((j+a_3)q, (j+a_3+1)q)$ is $\alpha+a_3(p-r_1)-p$ and so we must have $r_1<\alpha+a_3(p-r_1)-p<p$. Now returning to the first cluster of long segments, we see that the distance between $a_3q$ and the first multiple of $p$ in the interval $(a_3q, (a_3+1)q)$ is $a_3(p-r_1)$. Since the interval $(a_3q, (a_3+1)q)$ contains a long segment, we have that $0<a_3(p-r_1)<r_1$. But combining these inequalities leads to $p<\alpha$, a contradiction.

To prove property 7, let $x$ equal the number of clusters of long segments of length $a_3$ and $y$ equal the number of clusters of long segments of length $a_3+1$.  As in property 6, $r_1>1$ and $a_2=1$. From property 5, the number of short segments is $r_2+1$ and
hence $x+y=r_2$. Counting all the edges in the graph, we obtain
\begin{eqnarray*}
(r_2+1)a_1 + x a_3 (a_1+1) +y (a_3+1) (a_1+1) & = & q-1 \\
(a_1+1) + (r_2-y) a_3 (a_1+1) +y a_3 (a_1+1)+y(a_1+1) & = & q - r_2 a_1 \\
(a_1+1) + r_2 a_3 (a_1+1)+y(a_1+1) & = & a_1(r_1+r_2)+r_1-r_2 a_1 \\
(a_1+1) (1 + r_2 a_3 +y) & = &( a_1+1) r_1  \\
1 + r_1-r_3 +y & = & r_1  \\
y & = & r_3-1.
\end{eqnarray*}
Therefore, the number of clusters of long segments of length $a_3$ is $r_2-r_3+1$ and the number of clusters of length $a_3+1$ is $r_3-1$.
\hfill $\square$

\section{Relations among $\ell$-twisted Alexander polynomials}

Using Propositions~\ref{facts 1 about E} and \ref{facts 2 about E}, we can determine how $d_{p/q}(\lambda, t)$ and $e_{p/q}(\lambda, t)$ change when $p$ and $q$ are changed in certain ways. This will then allow us to compute the twisted Alexander polynomial for several different types of infinite families of 2-bridge knots. In order to do this, we need to investigate how $E(p,q)$ changes. 
It is useful to reformulate the contribution of each vertex in the graph $E(p,q)$ to $d$ and $e$ depending on the index of the vertex.  Examining the formulae for $d$ and $e$ given in (\ref{defn of d e}) again, we see that the $k$-th vertex  contributes $t^s \lambda^{\pm \lfloor\frac{k+1}{2}\rfloor}$ where $s$ depends on the $t$-level of the vertex as before and the plus or minus sign is determined as follows. If the vertex lies in the interior of a segment of slope $1$, then the choice is plus. If it lies in the interior of a segment of slope $-1$, then the choice is minus. If instead, it is a local maxima or minima, then $k=\sigma_1+\sigma_2+\dots+\sigma_i$, for some $i$, and the vertex contributes the term $t^{\sigma_1-\sigma_2+\dots+(-1)^{i-1}\sigma_i} \lambda^{(-1)^{(k+i+1)}\lfloor\frac{k+1}{2}\rfloor}$.

\begin{theorem}\label{increase clusters}  Suppose $1<p<q$, $pq$ odd, and $\gcd(p,q)=1$. Assume the Euclidean algorithm gives
$$\begin{array}{rclcl}
q&=&a_1 p+r_1, && 0<r_1<p\\
p&=&a_2 r_1+r_2,&& 0\le r_2<r_1\\
&\vdots&
\end{array}
$$
If $\ell$ is an odd prime dividing $q$, $p'=p+2 \ell jr_1$, and $q'=q+2 \ell ja_1 r_1$ for some positive integer $j$, then for the canonical forms from Theorem~\ref{twisted alex thm} we have
\vspace{ -.2 in }
\begin{enumerate}
\item $ \Delta_{K_{p'/q'}}(t) \equiv  \Delta_{K_{p/q}}(t) \ \mbox{\rm(mod $\ell$)}$ \ and   \  
$\tilde \Delta_{p'/q'}^\ell(t) \equiv  \tilde \Delta_{p/q}^\ell(t) \ \mbox{\rm(mod $\ell$).}$
\item If in addition, $\gcd(a_1, \ell)=1$, then $\displaystyle \frac{\tilde\Delta_{p'/q'}^{ \ell}(t)}{\Delta_{K_{p'/q'}}(t)} =\frac{\tilde\Delta_{p/q}^{ \ell}(t)}{\Delta_{K_{p/q}}(t)}$.
\end{enumerate}

\end{theorem}
{\bf Proof:} It suffices to prove the theorem for $j=1$. Notice that
$$\begin{array}{rclcl}
q'&=&a_1 p'+r_1, && 0<r_1<p\\
p'&=&(a_2+2 \ell) r_1+r_2,&& 0\le r_2<r_1.\\
\end{array}
$$
This implies that $1<p'<q'$, $p'q'$ odd, and $\gcd(p',q')=1$.

By Proposition~\ref{facts 2 about E}, the lengths of the short and long segments, respectively, of $E'=E(p',q')$ and $E=E(p,q)$ are the same. Secondly, the number of long segments is the same  in each graph. Moreover, the length of the initial cluster of short segments in $E'$ is $2\ell$ segments longer than that in $E$. Additionally, while the long segments may not be isolated in $E$, they certainly are in $E'$. We will now show that $E'$ is obtained from $E$ by increasing the length of every cluster of short segments by $2\ell$. In the case where there is a nonisolated long segment of $E$, and hence short segments are isolated, this will mean lengthening each cluster of short segments to a length of  $1+2\ell$ as well as  inserting clusters of short segments of length $2\ell$ between each adjacent pair of long segments.

Consider the $i$-th long segment in $E$. Suppose that it is preceded by $S$ short segments and $L=i-1$ long segments. This implies that
$$
(S+L)q<(Sa_1+L(a_1+1)+1)p < (Sa_1+L(a_1+1)+a_1+1)p<(S+L+1)q.
$$
These inequalities are equivalent to 
$$
(S'+L')q'<(S'a_1+L'(a_1+1)+1)p' < (S'a_1+L'(a_1+1)+a_1+1)p'<(S'+L'+1)q', 
$$
where $S'=S+(L+1)2\ell$ and $L'=L$.

This proves that each long segment in $E'$ is preceded by the correct number of long and short segments, namely, $L'=L$ and $S'=S+(L+1)2 \ell$, respectively. However, it  remains to show that the arrangement of short and long segments in $E'$ is obtained from that in $E$ as claimed. We prove this by induction on $i$. When $i=1$, we know by Proposition~\ref{facts 2 about E} that the first cluster of short segments is increased in length by $2\ell$. Assuming this is true up to the $(i-1)$-st long segment, and given that the $i$-th long segment is preceded by the correct number of long and short segments, we see that either the $i$-th cluster of short segments is increased in length by $2\ell$ or a new cluster of short segments of length $2\ell$ is inserted between the $(i-1)$-st and $i$-th long segments.

The first part of property 1 now follows since  lengthening  a single cluster of short segments by $2\ell$ will increase  two of the coefficients of the Alexander polynomial by $\ell$ and the $a_1-1$ remaining coefficients by $2\ell$. To prove the other  properties, we must examine how $d$ and $e$ are effected.

We now show that  $d$ and $e$ are the same for both $E$ and $E'$ if $a_1$ and $\ell$ are co-prime and are congruent modulo $\ell$ otherwise.  Since the length of each cluster of short segments is increased by $2\ell$ (and if the long segments are not isolated, new clusters of short segments of length $2\ell$  are introduced between adjacent long segments) we see that many new summands are introduced into both $d$ and $e$. However, each vertex in $E$ continues to contribute the same summand to either $d$ or $e$ of $E'$. To see this, consider the $k$-th vertex of $E$ which becomes the $k'$-th vertex of $E'$. Lengthening each cluster of short segments in $E$ by $2\ell$ means that $k'=k+2j\ell$ for some $j$.  Moreover, the vertex remains on the same $t$-level. If this vertex contributes $t^s \lambda^{\pm \lfloor \frac{k+1}{2}\rfloor}$ to $d$ or $e$ of $E$, then it will contribute $t^s \lambda^{\pm \lfloor \frac{k+1+2j\ell}{2}\rfloor}=t^s \lambda^{j\ell\pm \lfloor \frac{k+1}{2}\rfloor}=t^s \lambda^{\pm \lfloor \frac{k+1}{2}\rfloor}$ since $\lambda^{\ell}=1$.  Now we examine the contribution to $d$ or $e$ made by the insertion of a single set of $2\ell$ consecutive short segments. Consider a vertex of $E'$ in the  first of these $2\ell$ short segments. If this vertex is numbered $k$, then the reminding  $\ell-1$  vertices at the same $t$-level in this cluster are numbered $k+2a_1 i$ for $i=1, 2, \dots, \ell-1$.   These vertices all have $\lambda$ exponents of the same sign and together they contribute the following to either $d$ or $e$
\begin{eqnarray*}
&&t^s\left ( \lambda^{\pm \lfloor \frac{k+1}{2} \rfloor} +\lambda^{\pm(a_1+ \lfloor \frac{k+1}{2} \rfloor)} +\dots +\lambda^{\pm((\ell-1)a_1+ \lfloor \frac{k+1}{2} \rfloor)}               \right )\\
&=&t^s  \lambda^{\pm \lfloor \frac{k+1}{2} \rfloor} \left (1+ \lambda^{\pm a_1}+ ( \lambda^{\pm a_1})^2+\dots  +( \lambda^{\pm a_1})^{\ell-1}\right )\\
&=&\begin{cases}
0&\text{if $\gcd(a_1, \ell)=1$,}\\ 
\ell \, t^s  \lambda^{\pm \lfloor \frac{k+1}{2} \rfloor}& \text{otherwise.}
\end{cases}
\end{eqnarray*}
  \hfill $\square$

As an immediate consequence of Theorem~\ref{increase clusters} we have the following corollary.
\begin{corollary} Assume the hypotheses of Theorem~\ref{increase clusters} and additionally  that $\gcd(a_1, \ell)=1$.  If   $K_{p/q}$  satisfies  Conjecture~\ref{conjecture}, then so does $K_{p'/q'}$.
\end{corollary}

In Theorem~\ref{increase clusters} we found that $d$ and $e$ did not change when $p$ and $q$ were increased in a certain way. In the following result we keep $p$ fixed and increase $q$ in a certain way. This will now change $d$ and $e$ but in a such a way that will allow us to generate an infinite family of 2-bridge knots for which we can compute $d$ and $e$ recursively. This in turn will provide a recursive computation of $\ell$-twisted  Alexander polynomials.

\begin{proposition}\label{increase segments} Suppose $1<p<q$, $pq$  odd, $\gcd(p,q)=1$, and $\ell$ is an odd prime dividing $q$. Then for all $k>1$
\vspace{-.2 in}
\begin{enumerate}
\item \hspace{5 pt}$\Delta_{K_{p/( q+2 k \ell p)}}(t)=\alpha_k\, \Delta_{K_{p/( q+2\ell p)}}(t)+(1-\alpha_k) \Delta_{K_{p/q}}(t),$
\item $d_{p/(q+2 k \ell p)}(\lambda, t)=\alpha_k\, d_{p/(q+2\ell p)}(\lambda, t)+(1-\alpha_k) d_{p/q}(\lambda, t), \text{ and}$
\item \hspace{1 pt}$e_{p/(q+2 k \ell p)}(\lambda, t)=\alpha_k\, e_{p/(q+2\ell p)}(\lambda, t)+(1-\alpha_k) e_{p/q}(\lambda, t),$
\end{enumerate}
\vspace{-.2 in}
 where $\alpha_k=\frac{1-t^{2 k \ell}}{1-t^{2 \ell}}$.
\end{proposition}
\noindent{\bf Proof:} We focus on the second and third properties. The first property is proven in a similar way. Let $d=d_{p/ q}(\lambda, t), d'=d_{p/( q+2 \ell p)}(\lambda, t)$, and $d''=d_{p/(q+4\ell p)}(\lambda, t)$. Define $e, e'$, and $e''$ analogously. We will show that
\begin{equation*}d''-d'=t^{2\ell}(d'-d) \quad \text{ and } \quad e''-e'=t^{2\ell}(e'-e).\end{equation*}
This will imply that
$$d''=(1+t^{2 \ell})d'-t^{2 \ell}d \quad \text{ and } \quad  e''=(1+t^{2 \ell})e'-t^{2 \ell}e.$$
Standard techniques of linear algebra can then be used to derive the recursion formula stated in the lemma. 

Notice that $d$ and $e$ can be uniquely recovered from the sum $d+e$ by separating the terms of even and odd powers of $t$. Hence it suffices to prove that 
\begin{equation}\label{d+e difference} 
(d''+e'')-(d'+e')=t^{2\ell}[(d'+e')-(d+e)].
\end{equation}

Let $\sigma=\sigma(p, q)=\{\sigma_1, \sigma_2, \dots, \sigma_p\}$, $\sigma'=\sigma(p, q+2 \ell p)=\{\sigma'_1, \sigma'_2, \dots, \sigma'_p\}$, $\sigma''=\sigma(p, q+4\ell p)=\{\sigma''_1, \sigma''_2, \dots, \sigma''_p\}$, and  let 
$E, E'$, and $E''$ be the associated epsilon graphs. It follows from property 5 of Proposition~\ref{facts 1 about E} that  $\sigma'_i=\sigma_i+2 \ell$ and $\sigma''_i=\sigma'_i+2 \ell$. Hence $E'$ is obtained from $E$ by increasing each segment by $2\ell$ edges. Similarly,  $E''$ is obtained from $E'$ by increasing each segment by $2\ell$ edges. Superimposing the graphs, we see that  the local minima of all three graphs coincide. 

The $k$-th vertex of  $E$ contributes a term of the form $t^s \lambda^{\pm \lfloor \frac{k+1}{2}\rfloor}$ to  $d+e$ as follows. For the moment, we focus only on the exponent of $\lambda$. To describe the sign of the exponent of $\lambda$ better, we consider three cases: the vertex is in the interior of a segment with slope 1, the vertex is in the interior of a segment of slope $-1$, or  the vertex is a local extrema. In the first case, the exponent of $\lambda$ is $\lfloor \frac{k+1}{2}\rfloor$ and in the second case it is $-\lfloor \frac{k+1}{2}\rfloor$. If the vertex is a local extrema, then 
we must have that $k=\sigma_1+\sigma_2+\dots+\sigma_i$ for some $i$. The exponent of $\lambda$ associated to this term is now $(-1)^{k+i+1} { \lfloor \frac{k+1}{2}\rfloor}$. Notice that $k$ is the total number of edges that precede the vertex while $i$ is the total number of segments that precede the vertex.

Suppose $v$ is a vertex that is a local minima of all three epsilon graphs. If $v$ is the $k$-th vertex of $E$, then $k=\sigma_1+\sigma_2+\dots+\sigma_{2i}$ for some $i$. If $v$ is the $k'$-th and $k''$-th vertex of $E'$ and $E''$, respectively, then $k\equiv k' \equiv k''$ (mod $2\ell$) since $\sigma''_j \equiv \sigma'_j \equiv \sigma_j$ (mod $2\ell$). Hence $v$ contributes the same term to $d+e$, $d'+e'$, and $d''+e''$ because $\lambda^{\ell}=1$. Thus these terms cancel in the differences on both sides of (\ref{d+e difference}).

Consider now a local minima shared by the three epsilon graphs and the three segments of slope 1 that start at this minima and belong to the three graphs, respectively. Suppose $v$ and $v'$ are two vertices in the interior of the segments of $E$ and $E'$, respectively, that lie on the same $t$-level, say, $t^s$. Suppose $v$ is the $k$-th vertex of $E$ and $v'$ is the $k'$-th vertex of $E'$. Again, since $\sigma'_j \equiv \sigma_j$ (mod $2\ell$), we see that  $k'\equiv k$ (mod $2\ell$). The vertex $v$ contributes the term $t^s \lambda^{ \lfloor \frac{k+1}{2}\rfloor}$ to $d+e$ while the vertex $v'$ contributes the term $t^s \lambda^{ \lfloor \frac{k'+1}{2}\rfloor}$ to $d'+e'$. But, because  $k'\equiv k$ (mod $2\ell$) and $\lambda^{\ell}=1$, we see that the two contributions are the same. Thus they will cancel in the difference given in (\ref{d+e difference}). The same is true of any two vertices at the same $t$-level but in  the interior of segments of $E''$ and $E'$. Moreover, terms of this kind that lie at the same $t$-level on the interiors of segments of slope $-1$ that share the same local minima also cancel. 

Next, we consider a local maxima of $E$. For such a vertex $k=\sigma_1+\sigma_2+\dots \sigma_{2i+1}$. The corresponding maxima of $E'$ is the vertex numbered $k'=\sigma'_1+\sigma'_2+\dots \sigma'_{2i+1}$.  If this is not the last vertex of $E$, then there exist two vertices of $E'$  at the same $t$-level, namely the vertices numbered $k'-2\ell$ and $k'+2\ell$. The first contributes 
$t^s \lambda^{ \lfloor \frac{k'-2\ell+1}{2} \rfloor}$ to $d'+e'$ 
and the second contributes $t^s \lambda^{- \lfloor \frac{k'+2\ell+1}{2} \rfloor}$  to $d'+e'$. These exponents of $\lambda$ differ only in sign. Now, returning to the local maxima of $E$, we see that it contributes $t^s \lambda^{\pm \lfloor \frac{k+1}{2} \rfloor}$ and hence cancels with one of the two terms from $E'$. If in fact, the local maxima of $E$ is the last vertex of $E$, then it contributes $t^s \lambda^{\lfloor \frac{k+1}{2} \rfloor}$ and again cancels with the term of $E'$. In conclusion we see that every term of $d+e$ cancels with some term of $d'+e'$.

Finally, it remains to see that each term of $E''$ that does not cancel with a term of $E'$ is equal to $t^{2\ell}$ times a corresponding term of $E'$ that does not cancel with a term of $E$. Consider for example  $k''=\sigma''_1+\sigma''_2+\dots \sigma''_{2i+1}$ and the vertex $v''$ of $E''$ numbered $k''-j$ for some $0\le j<2\ell$. Let $v'$ be the vertex of $E'$ numbered $k'-j$ where   $k'=\sigma'_1+\sigma'_2+\dots \sigma'_{2i+1}$. If the local maxima of $E''$ at position $k''$ has $t$-level $s$, then $v''$ contributes the term $t^{s-j}\lambda^{\lfloor \frac{k''-j+1}{2} \rfloor}$ to $(d''+e'')-(d'+e')$ while the vertex $v'$ of $E'$ contributes $t^{s-2\ell-j}\lambda^{\lfloor \frac{k'-j+1}{2} \rfloor}$ to $(d'+e')-(d+e)$. Similar pairings of vertices exist for the segments of slope $-1$.
\hfill $\square$

\section{Special families of $\ell$-twisted Alexander polynomials}

In this section we apply Theorems~\ref{twisted alex thm} and \ref{increase clusters} as well as Proposition~\ref{increase segments} to compute the $\ell$-twisted Alexander polynomials for several infinite families of 2-bridge knots. Obtaining a specific formula for the $\ell$-twisted Alexander polynomial allows us to verify Conjecture~\ref{conjecture} for these families. We begin with 2-bridge torus knots. 

\subsection{2-bridge torus knots}

If $p=1$, the 2-bridge knot $K_{1/q}$ is the  $(2, q)$-torus knot. Conjecture~\ref{conjecture} was proven for these knots in the case where $q$ is prime by Hirasawa and Murasugi  \cite{HM:2009}. However,  they do this without producing an explicit formula for the twisted Alexander polynomial even though   they provide a likely candidate. (See Remark 5.3 of  \cite{HM:2009}.) The following theorem extends their result to any odd integer $q$ and moreover confirms their conjectured formula in the case when $q$ is prime.

\begin{theorem}\label{torus knots} If $\ell$ is an odd prime dividing $q$, then the $\ell$-twisted Alexander polynomial of the 2-bridge, torus knot $K_{1/q}$ is,  
$$\tilde \Delta_{1/ q}^{\ell}(t) \doteq \frac{\Delta_{K_{1/q}}(t)}{t-1} \left( \frac{(1+t)(1+t^q)^\frac{\ell-1}{2}}{1+t^{\ell}} \right)\left( \frac{(1-t)(1-t^q)^\frac{\ell-1}{2}}{1-t^{\ell}} \right).$$
\end{theorem}
\vspace{-.3 in}
\noindent{\bf Proof:}
Let $\ell$ be an odd prime dividing $q$. Since $\epsilon_i=1$ for all $i$,  we have from (\ref{defn of d e}) that
$$d_{1/q}(\lambda,t) = \sum_{k=1}^\frac{q-1}{2} t^{2k-1} \lambda^{k} = \frac{t \lambda(t^{q-1}\lambda^\frac{q-1}{2}-1)}{t^2 \lambda-1} \quad 
\text{and} \quad  
e_{1/q}(\lambda,t) = 1+\sum_{k=1}^\frac{q-1}{2} t^{2k} \lambda^k = \frac{t^{q+1} \lambda^\frac{q+1}{2}-1}{t^2 \lambda-1}.$$
From this and the identity $\lambda^\frac{q+1}{2} = \lambda^\frac{-q+1}{2}$ we see that
$$
D_{1/q}(\lambda, t) = -\frac{(t^{2q}-1)(t^2-1)}{(t^2 -\lambda)(t^2 -\lambda^{-1})}.
$$
This leads to 
\begin{eqnarray*}
\mathlarger{\mathlarger{\prod}}_{i=1}^{\frac{\ell-1}{2}}D_{1/q}(\lambda^i, t)&=&
\mathlarger{\mathlarger{\prod}}_{i=1}^{\frac{\ell-1}{2}} \frac{-(t^{2q}-1)(t^2-1)}{(t^2- \lambda^i)(t^2- \lambda^{-i})}\\
&=&\mathlarger{\mathlarger{\prod}}_{i=1}^{\frac{\ell-1}{2}} \frac{-(t^{2q}-1)(t^2-1)}{(t^2- \lambda^i)(t^2- \lambda^{\ell-i})}\\
&=&\frac{(-1)^\frac{\ell-1}{2} (t^{2q}-1)^\frac{\ell-1}{2}(t^2-1)^\frac{\ell-1}{2}}{\prod_{i=1}^{\ell-1}(t^2 -\lambda^i)}.
\end{eqnarray*}

If $\omega$ is a primitive $\ell$-th root of unity, then
\begin{equation}\label{primitive root identity}
\frac{x^{\ell}-1}{x-1}=x^{\ell-1} + \dots +x^2+x+1=(x-\omega)(x-\omega^2)\dots (x-\omega^{\ell-1}).
\end{equation}
Setting $x=t^2$ and $\omega=\lambda$ leads to

$$\prod_{i=1}^{\ell-1} (t^2-\lambda^i)=\frac{t^{2 \ell}-1}{t^2-1}.$$
Combining these results with Theorem~\ref{twisted alex thm} completes the proof.
\hfill $\square$

\begin{corollary}
Conjecture~\ref{conjecture} is true for all 2-bridge, torus knots.
\end{corollary}
\noindent{\bf Proof:}
Letting
$$f(t) = \frac{(1+t)(1+t^q)^\frac{\ell-1}{2}}{1+t^{\ell}}, $$
the first part of the conjecture is clear from Theorem~\ref{torus knots}.
Furthermore, since $\ell$ is prime, it follows that $(1+t)^{\ell} \equiv 1+t^{\ell} \ (\mbox{mod } \ell)$. Hence,
\begin{eqnarray*}
f(t) & \equiv & \frac{(1+t)(1+t^q)^\frac{\ell-1}{2}}{(1+t)^{\ell}} \\
& \equiv & \frac{(1+t^q)^\frac{\ell-1}{2}}{(1+t)^\frac{\ell-1}{2}(1+t)^\frac{\ell-1}{2}} \\
& \equiv & \left( \frac{(1+t^q)}{(1+t)(1+t)} \right)^\frac{\ell-1}{2} \ \ (\mbox{mod } \ell).
\end{eqnarray*}
Using (\ref{alex poly formula}), we obtain
$$\Delta_{K_{1/q}}(t)=1-t+t^2-\dots+t^{q-1}=\frac{1+t^q}{1+t}.$$
Therefore,  these choices of $\Delta_{K_{1/q}}(t)$ and $f(t)$ also satisfy the second part of the conjecture. 
\hfill$\square$

\subsection{2-bridge knots with genus one}

In this section we compute the $\ell$-twisted Alexander polynomial  for all 2-bridge knots with genus one. The 2-bridge knot $K_{p/q}$ has genus one if and only if 
$$\frac{p}{q} = \frac{4rs-2s \pm1}{4 rs\pm1},$$ for some natural  numbers $r$ and $s$.
(See  Proposition~12.25 of  \cite{BZ:2003}.) We begin by determining the epsilon  graph for these fractions. The following lemma follows easily from the facts in Propositions~\ref{facts 1 about E} and \ref{facts 2 about E}.

\begin{lemma}
\label{epsilons for genus 1}Suppose $r$ and $s$ are natural numbers. Then
\vspace{-.2 in}
\begin{enumerate}
\item $\sigma(4rs-2s-1, 4rs-1)=\{[1]_{2r-2}, 2, [1]_{2r-2}, 2, [1]_{2r-2}, 2, \dots, [1]_{2r-2}, 2, [1]_{2r-2}\}$ and
\item 
$\sigma(4rs-2s+1, 4rs+1)=\{[1]_{2r-1}, 2, [1]_{2r-2}, 2,  [1]_{2r-2}, 2, \dots,  [1]_{2r-2}, 2,  [1]_{2r-1}\}$.
 \end{enumerate}
 \vspace{-.2 in}

 where $[n]_k=\underbrace{n,n, \dots, n}_k$ and the number of 2's in each case is $2s-1$.

\end{lemma}
 
Using Theorem~\ref{twisted alex thm} and Lemma~\ref{epsilons for genus 1} we are now prepared to find the $\ell$-twisted Alexander polynomial of 2-bridge knots with genus one.

\begin{theorem}\label{genus 1 knots} If $K_{p/q}$ is a 2-bridge knot with genus one and $\ell$ is an odd prime dividing $q$, then
\begin{enumerate}
\item $\displaystyle \tilde\Delta_{p,q}^{\ell}(t)  \doteq  \frac{\Delta_{K_{p/q}}(t)}{t-1}(1+t)^\frac{\ell-1}{2}(1-t)^\frac{\ell-1}{2},  \hspace{30 pt} \text{ if } \   \frac{p}{q} = \frac{4rs-2s -1}{4 rs-1} \text{ and}$\\
\item $\displaystyle  \tilde\Delta_{p,q}^{\ell}(t)  \doteq     \frac{\Delta_{K_{p/q}}(t)}{t-1} \left (\frac{1+t}{t}\right )^\frac{\ell-1}{2}\left (\frac{1-t}{-t}\right )^\frac{\ell-1}{2}
,  \ \text{ if } \   \frac{p}{q} = \frac{4rs-2s +1}{4 rs+1}$.
\end{enumerate}
\end{theorem}

{\bf Proof:} Suppose first that $\frac{p}{q}=\frac{4rs-2s-1}{4rs-1}$. Then from Lemma~\ref{epsilons for genus 1} we have 
$$\sigma(p, q)=\{[1]_{2r-2}, 2, [1]_{2r-2}, 2, [1]_{2r-2}, 2, \dots, [1]_{2r-2}, 2, [1]_{2r-2}\}$$
where the number of isolated 2's is $2s-1$.  Now label the epsilon graph as described at the end of Section~\ref{twisted alex poly}.  Inspection of the vertices on the middle horizontal level of the epsilon graph gives
\begin{eqnarray*}
d_{p/q}(\lambda,t)&=&t \big[ (\lambda^{-1}+ \dots + \lambda^{-r+1} )+(\lambda^{r}+\dots+\lambda^{2r-1})+ \\
&&\phantom{t \big[} (\lambda^{-2r}+\dots+\lambda^{-3r+1}) + (\lambda^{3r}+\dots+\lambda^{4r-1}) + \dots \\
&& \phantom{t \big[} (\lambda^{-2rs+2r}+\dots+\lambda^{-2rs+r+1})+(\lambda^{2rs-r}+\dots+\lambda^{2rs-1})].
\end{eqnarray*}
Rewriting this expression, multiplying by $(1-\lambda)/(1-\lambda)$, and then simplifying gives
\begin{eqnarray*}
d_{p/q}(\lambda,t)&=&t \big[ -1 + (1+ \dots + \lambda^{-r+1} )+(\lambda^{-2r}+\dots+\lambda^{-3r+1}) + \dots + \\
&& \phantom{t \big[} (\lambda^{-2rs+2r}+\dots+\lambda^{-2rs+r+1}) + (\lambda^{r}+\dots+\lambda^{2r-1})+ \\
&&\phantom{t \big[} (\lambda^{3r}+\dots+\lambda^{4r-1}) + \dots +(\lambda^{2rs-r}+\dots+\lambda^{2rs-1})]\\
&=&\frac{t}{1-\lambda} \bigg[-(1-\lambda)+ (\lambda^{-r+1}-\lambda) +(\lambda^{-3r+1}-\lambda^{-2r+1})+ \dots + \\
&&\phantom{\frac{t}{1-\lambda} \bigg[}(\lambda^{-2rs+r+1}-\lambda^{-2rs+2r+1})+(\lambda^{r}-\lambda^{2r})+(\lambda^{3r}-\lambda^{4r})+\dots+(\lambda^{2rs-r}-\lambda^{2rs})\bigg].
\end{eqnarray*}
Next we multiply by $(1+\lambda^r)/(1+\lambda^r)$ and use the fact that $\lambda^{4rs-1}=1$ to obtain
\begin{eqnarray*}
d_{p/q}(\lambda,t)&=&\frac{t}{(1-\lambda)(1+\lambda^r)} \bigg[-(1+\lambda^r)(1-\lambda)+ (\lambda^{-2rs+r+1}-\lambda^{r+1})+(\lambda^r-\lambda^{2rs+r})\bigg]\\
&=&\frac{t}{(1-\lambda)(1+\lambda^{r})} \bigg[-(1+\lambda^r)(1-\lambda)+ \lambda^{2rs+r}-\lambda^{r+1}+\lambda^r-\lambda^{2rs+r} \bigg]\\
&=&\frac{t}{(1-\lambda)(1+\lambda^{r})} \bigg[-(1+\lambda^r)(1-\lambda)+\lambda^r(1-\lambda) \bigg]\\
&=&\frac{-t}{1+\lambda^{r}}.
\end{eqnarray*}
Returning to the labeled epsilon graph, we find that $e_{p/q}(\lambda,t)=\alpha+\beta t^2$, where $\alpha$ and $\beta$ are Laurent polynomials in $\lambda$. Furthermore, using the rotational symmetry of the epsilon graph, we see that $\alpha \lambda^{2rs-1} =\beta$, and so, $\alpha = \beta \lambda^{2rs}$. Reading $\beta$ off from the graph, we have
\begin{eqnarray*}
\beta&=&(\lambda^r+\lambda^{r+1}+\dots+\lambda^{2r-1})+(\lambda^{3r}+\lambda^{3r+1}+\dots+\lambda^{4r-1})+\dots+\\
&&(\lambda^{2rs-r}+\lambda^{2rs-r+1}+\dots+\lambda^{2rs-1})\\
&=&\frac{1}{1-\lambda}[(\lambda^r-\lambda^{2r})+(\lambda^{3r}-\lambda^{4r})+\dots+(\lambda^{2rs-r}-\lambda^{2rs})]\\
&=&\frac{1 - \lambda^{2rs}}{(1-\lambda)(1+\lambda^{-r})}.
\end{eqnarray*}
Hence
$$e_{p/q}(\lambda,t)=\frac{(\lambda^{2rs}+t^2)(1-\lambda^{2rs})}{(1-\lambda)(1+\lambda^{-r})}.$$

Using these formulas, we may now compute $D_{p/q}(\lambda, t)$.
\begin{eqnarray*}
D_{p/q}(\lambda, t)&=&\frac{t^2}{(1+\lambda^r)(1+\lambda^{-r})}-\frac{(\lambda^{2rs}+t^2)(\lambda^{-2rs}+t^2)(1-\lambda^{2rs})(1-\lambda^{-2rs})}{(1-\lambda)(1-\lambda^{-1})(1+\lambda^r)(1+\lambda^{-r})}\\
&=&\frac{(1-\lambda^{4rs})(1-\lambda^{-4rs})t^2-(1+(\lambda^{2rs}+\lambda^{-2sr})t^2+t^4)(1-\lambda^{2sr})(1-\lambda^{-2sr})}{(1-\lambda)(1-\lambda^{-1})(1+\lambda^r)(1+\lambda^{-r})}\\
&=&-\frac{(1-\lambda^{2sr})(1-\lambda^{-2sr})(t+1)^2(t-1)^2}{(1-\lambda)(1-\lambda^{-1})(1+\lambda^r)(1+\lambda^{-r})}.
\end{eqnarray*}
Here we have used the fact that $\lambda^{4rs}=\lambda$. Now
\begin{eqnarray*}\label{product of determinants}
\mathlarger{\mathlarger{\prod}}_{i=1}^{\frac{\ell-1}{2}}D_{p/q}(\lambda^i, t)&=&(-1)^\frac{\ell-1}{2}(1+t)^{\ell-1}(1-t)^{\ell-1}\mathlarger{\mathlarger{\prod}}_{i=1}^{\frac{\ell-1}{2}} 
\frac{(1-\lambda^{2sri})(1-\lambda^{-2sri})}{(1-\lambda^i)(1-\lambda^{-i})(1+\lambda^{ri})(1-\lambda^{-ri})}\\
&=&(-1)^\frac{\ell-1}{2}(1+t)^{\ell-1}(1-t)^{\ell-1}\mathlarger{\mathlarger{\prod}}_{i=1}^{\ell-1}\frac{1-\lambda^{2sri}}{(1-\lambda^i)(1+\lambda^{ri})}.
\end{eqnarray*}
Setting $x=1$ and $x=-1$ in (\ref{primitive root identity}), we obtain the following two identities:
$$\prod_{i=1}^{\ell-1}(1-\omega^i)=\ell \qquad \text{and} \qquad \prod_{i=1}^{\ell-1}(1+\omega^i)=1.$$
Since both $r$ and $2rs$ are relatively prime to $\ell$, it follows that both $\lambda^r$ and $\lambda^{2rs}$ are also primitive. Hence we obtain
\begin{equation*}
\mathlarger{\mathlarger{\prod}}_{i=1}^{\ell-1}\frac{1-\lambda^{2sri}}{(1-\lambda^i)(1+\lambda^{ri})}=1.
\end{equation*}

Substituting these results into (\ref{twisted poly formula}) now gives 
\begin{eqnarray*}
\tilde \Delta_{p/q}^{\ell}(t) & \doteq & \frac{\Delta_{K_{p/q}}(t)}{(t+1)^\frac{\ell-1}{2}(t-1)^\frac{\ell+1}{2}} (-1)^\frac{\ell-1}{2}(1+t)^{\ell-1}(1-t)^{\ell-1} \\
& \doteq & \frac{\Delta_{K_{p/q}}(t)}{t-1} (t+1)^\frac{\ell-1}{2}(t-1)^\frac{\ell-1}{2}.
\end{eqnarray*}

Using a similar analysis for the case $\frac{p}{q}=\frac{4rs-2s+1}{4rs+1}$ we find that 
$$e_{p/q}(\lambda,t) = \frac{1}{1+\lambda^r} \ \ \ \  \mbox{and} \ \ \ \ d_{p/q}(\lambda,t) =\frac{(t^{-1}+\lambda^{2rs}t) (1-\lambda^{2rs})\lambda}{(1-\lambda)(1+\lambda^{-r})}.$$
From this we have
$$D_{p/q}(\lambda, t)=\frac{(1-\lambda^{2sr})(1-\lambda^{-2sr})(t+1)^2(t-1)^2}{t^2(1-\lambda)(1-\lambda^{-1})(1+\lambda^r)(1+\lambda^{-r})}$$
which is the determinant in the previous case divided by  $-t^2$.  This gives the second formula.   \hfill $\square$

We conclude by demonstrating that the conjecture of Hirasawa and Murasugi is indeed true for 2-bridge knots with genus one. 
\begin{corollary} Conjecture~\ref{conjecture} is true for all 2-bridge knots with genus one.
\end{corollary}
\noindent{\bf Proof:}
Letting $f(t)=(1+t)^\frac{\ell-1}{2}$, the first part of the conjecture is clear from Theorem~\ref{genus 1 knots}. Thus, it remains to show that 
$$f(t)\equiv \left (\frac{\Delta_{K_{p/q}}(t)}{1+t}\right)^\frac{\ell-1}{2} (\text{mod } \ell).$$
If $\frac{p}{q}=\frac{4rs-2s-1}{4rs-1}$, then it is easy to verify from the epsilon graph that 
$$\Delta_{K_{p/q}}(t)=rs-(2rs-1)t+rs t^2.$$
Since $\ell$ divides $4rs-1$, we have $4rs \equiv 1 \ (\text{mod } \ell)$ and hence 
\begin{eqnarray*}
\Delta_{K_{p/q}}(t)&\equiv&rs+2rst+rst^2  \\
&\equiv&rs(1+t)^2  \ (\text{mod }\ell).
\end{eqnarray*}
Also, by Fermat's Little Theorem, we have $2^{\ell-1}\equiv 1\ (\text{mod } \ell)$; hence
\begin{eqnarray*}
1&\equiv &(4rs)^\frac{\ell-1}{2} \\
&\equiv & 2^{\ell-1}(rs)^\frac{\ell-1}{2} \\
&\equiv & (rs)^\frac{\ell-1}{2} \ (\text{mod } \ell).
\end{eqnarray*}
Therefore,
\begin{eqnarray*}
\left(\frac{\Delta_K(t)}{1+t}\right)^\frac{\ell-1}{2} &\equiv &\left(\frac{rs(1+t)^2}{1+t}\right)^\frac{\ell-1}{2}\\
&\equiv& (rs)^\frac{\ell-1}{2}(1+t)^\frac{\ell-1}{2} \\
&\equiv&f(t) \ (\text{mod } \ell).
\end{eqnarray*}

In the case of $\frac{p}{q}=\frac{4rs-2s+1}{4rs+1}$, we have   
$\Delta_{K_{p/q}}(t)=-rst^{-1}+(2rs+1)-rs t$, $f(t)=t^\frac{1-l}{2}(t+1)^\frac{l-1}{2}$, 
and the proof is similar.
\hfill $\square$

\subsection{Recursive Families}

If we apply the change in $p$ and $q$ given in Theorem~\ref{increase clusters} to a genus-one 2-bridge knot, we produce another genus-one 2-bridge knot. Hence we cannot use Theorem~\ref{increase clusters} to extend the results of Theorem~\ref{genus 1 knots}. In this sub-section, we show how to use Theorem~\ref{increase clusters} and Proposition~\ref{increase segments}   to derive formulae for $\ell$-twisted Alexander polynomials of various bi-infinite families of 2-bridge knots. Using these formulae we can confirm Conjecture~\ref{conjecture}. From a practical point of view, we can do this for any 2-bridge knot with  fixed values of $p$ and $\ell$. We begin with the case of $p=5$ and $\ell=3$.

Given any 2-bridge knot $K_{p/q}$ with $p=5$ and $q$ divisible by $\ell=3$, if we reduce $q$ by multiples of $2\ell p=30$ we arrive at one of the following {\it root} fractions: $5/9,\ 5/21,\ 5/27$, or $5/33$. Consider first the fraction $5/9$. Computing $d$ and $e$ for $5/9$ and $5/39$, and using the fact that $\lambda^3=1$,  we obtain
$$\begin{array}{rclrcl}
d_{5/9}(\lambda, t)&=&-t^{-1}-t \lambda,  \quad &e_{5/9}(\lambda, t)&=&-\lambda\\
d_{5/39}(\lambda, t)&=&-t^{-1}-t \lambda +t^3(1+\lambda)-t^5-t^7 \lambda,  \quad &e_{5/39}(\lambda, t)&=&-\lambda-t^2(1+\lambda)-t^4+t^6\lambda.\\
\end{array}$$

Propostion~\ref{increase segments} now tells us that
\begin{eqnarray*}
d_{5/(9+30k)}(\lambda, t)&=&\alpha_k \ d_{5/39}(\lambda, t)+(1-\alpha_k)\, d_{5/9}(\lambda, t)\\
e_{5/(9+30k)}(\lambda, t)&=&\alpha_k\  e_{5/39}(\lambda, t)+(1-\alpha_k)\, e_{5/9}(\lambda, t)
\end{eqnarray*}
where $\alpha_k=\frac{1-t^{6k}}{1-t^6}$. If we now simplify the determinant $D_{5/(9+30k)}(\lambda, t)$
we obtain
$$
D_{5/(9+30k)}(\lambda, t)=\frac{(t-1)^2 (t+1)^2 \left(t^6 \alpha _k-t^3 \alpha _k+1\right)\left(t^6 \alpha _k+t^3 \alpha _k+1\right)}{t^2}
$$
Thus the twisted Alexander polynomial is
$$\tilde \Delta_{5/(9+30k)}^3(t)\doteq \frac{\Delta_{5/(9+30k)}(t)}{t-1}t^4(t-1)(t+1) \left(t^6 \alpha _k-t^3 \alpha _k+1\right)
   \left(t^6 \alpha _k+t^3 \alpha _k+1\right).$$
 If we let $f(t)=t^2 (t+1)(t^6 \alpha_k- t^3 \alpha_k+1)$, then the polynomial has the form given in Conjecture~\ref{conjecture}. To verify the rest of the conjecture we need to show that 
$$f(t)\equiv \frac{\Delta_{5/(9+30 k)}(t)}{1+t} \quad \text{(mod 3)}.$$

Computing the Alexander polynomials of $K_{5/9}$ and $K_{5/39}$ we obtain
\begin{eqnarray*}
\Delta_{5/9}(t)&=&-2t^{-1}+5-2t,\\
 \Delta_{5/39}(t)&=&-2t^{-1}+5-5t+5t^2-5t^3+5t^4-5t^5+5t^6-2t^7.
 \end{eqnarray*}
Hence by Proposition~\ref{increase segments},
 $$\Delta_{5/(9+30k)}(t)\doteq t^3\left (\alpha_k (-2t^{-1}+5-5t+5t^2-5t^3+5t^4-5t^5+5t^6-2t^7)+(1-\alpha_k)(-2t^{-1}+5-2t)\right ).$$
 It is now straightforward to check that this polynomial, divided by $1+t$, is congruent modulo 3 to $f(t)$.  Since $9=1\cdot 5+4$ and $\gcd(1, 3)=1$, we can combine this result with Theorem~\ref{increase clusters} to prove the following Theorem.

 \begin{theorem} Let $k$ and $j$ be non-negative integers. Then
$$\tilde \Delta_{(5+24 j)/(9+30k+144 jk+24 j)}^3(t)\doteq \frac{\Delta_{K}(t)}{t-1}t^4(t-1)(t+1) \left(t^6 \alpha _k-t^3 \alpha _k+1\right)
   \left(t^6 \alpha _k+t^3 \alpha _k+1\right),$$
   where $\alpha_k=\frac{1-t^{6k}}{1-t^6}$.
   Moreover, these knots satisfy Conjecture~\ref{conjecture} with  $f(t)=t^2 (1+t)(1- t^3 \alpha_k+t^6 \alpha_k)$.
\end{theorem}
 
Continuing with $p=5$ and $\ell=3$, there are three more families to consider, namely those with root fractions $5/21$, $5/27$, and $5/33$. However, in the case of $5/33$ we have $33=6\cdot 5+3$ and $\gcd(6,3)\ne1$. Hence we may apply Propostion~\ref{increase segments} but not Theorem~\ref{increase clusters} to the root fraction $5/33$. For the root fractions $5/21$ and $5/27$, the 3-twisted Alexander polynomial of all the knots in the bi-infinite families indexed by $j$ and $k$  can be computed in a way completely similar to the case of $5/9$ and used to verify Conjecture~\ref{conjecture}. For the root fraction $5/33$ we obtain a similar result but must assume that $j=0$. We omit the details and summarize our results in Theorem~\ref{additional results} and the Appendix. We include in the Appendix not only the remaining three families with $p=5$ and $\ell=3$, but a sampling of several more families.

\begin{theorem}\label{additional results} For each  root fraction $p/q$ listed in the Appendix,  the $\ell$-twisted Alexander polynomial of $K_{p'/q'}$, where 
$$\frac{p'}{q'}=\frac{p+2 \ell j r}{q+2\ell(kp+jar+2 \ell j k r)},$$ 
for  any non-negative integers $k$ and $j$, is of the form
$$\tilde \Delta_{p'/q'}^\ell(t)\doteq \frac{\Delta_{K_{p'/q'}}(t)}{t-1}f(t)f(-t)$$
where $f(t)$ is given in the Appendix and $\alpha_k=\frac{t^{2 \ell k}-1}{t^{2 \ell}-1}$. Moreover
$$f(t)\equiv \left ( \frac{\Delta_{K_{p'/q'}}(t)}{t+1}\right )^\frac{\ell-1}{2} \quad \text{{\rm (mod} $\ell$\rm)}.$$
If the root fraction appears in bold-face, then we must assume that $j=0$.
\end{theorem}


\pagebreak
{\bf Appendix}

The following table gives the Laurent polynomial $f(t)$ described in Theorem~\ref{additional results} for given root faction $p/q$ and odd prime $\ell$.

 \renewcommand{\arraystretch}{.5}
 \begin{center}$$\begin{array}{ccl}
p/q &\ell&f(t)\\ 
 \\ 
\hline 
\\ 
5/9&3&t^{-1}\left(t+1\right)\left(t^{6}\alpha_k-t^{3}\alpha_k+1\right)\\ 
\\ 
5/21&3&t\left(t+1\right)\left(t^{6}\alpha_k-t^{3}\alpha_k+1\right)\\ 
\\ 
5/27&3&-\left(t+1\right)\left(t^{10}\alpha_k-t^{9}\alpha_k-t^{8}\alpha_k-t^{7}\alpha_k+t^{6}\alpha_k+t^{5}\alpha_k+t^{4}-t^{3}-t^{2}-t+1\right)\\ 
\\ 
{\bf 5/33}&3&t^{-2}\left(t+1\right)\left(t^{12}\alpha_k+t^{11}\alpha_k-t^{10}\alpha_k-t^{9}\alpha_k-t^{8}\alpha_k+t^{7}\alpha_k+t^{6}+t^{5}-t^{4}-t^{3}-t^{2}+t+1\right)\\ 
 \\ 
\hline 
\\ 
7/9&3&t^{-1}\left(t+1\right)\left(t^{6}\alpha_k+t^{5}\alpha_k-t^{4}\alpha_k-t^{3}\alpha_k-t^{2}\alpha_k+t\alpha_k+1\right)\\ 
\\ 
7/15&3&\left(t+1\right)\left(t^{6}\alpha_k-t^{3}\alpha_k+1\right)\\ 
\\ 
{\bf 7/27}&3&\left(t+1\right)\left(t^{6}\alpha_k-t^{3}\alpha_k+1\right)\\ 
\\ 
7/33&3&-t^{-1}\left(t+1\right)\left(t^{10}\alpha_k-t^{9}\alpha_k-t^{8}\alpha_k-t^{7}\alpha_k+t^{6}\alpha_k+t^{5}\alpha_k+t^{4}-t^{3}-t^{2}-t+1\right)\\ 
\\ 
7/39&3&-t^{-3}\left(t+1\right)\left(t^{12}\alpha_k+t^{11}\alpha_k-t^{10}\alpha_k-3t^{9}\alpha_k-t^{8}\alpha_k+t^{7}\alpha_k+2t^{6}\alpha_k+t^{6}+t^{5}-t^{4}\right. \\ 
 \\ 
 &&\left. \hspace{.1 in}-3t^{3}-t^{2}+t+1\right)\\ 
\\ 
{\bf 7/45}&3&\left(t+1\right)\left(2t^{12}\alpha_k+t^{11}\alpha_k-t^{10}\alpha_k-3t^{9}\alpha_k-t^{8}\alpha_k+t^{7}\alpha_k+t^{6}\alpha_k+2t^{6}+t^{5}-t^{4}-3t^{3}-t^{2}\right. \\ 
 \\ 
 &&\left. \hspace{.1 in}+t+2\right)\\ 
 \\ 
\hline 
\\ 
3/5&5&t^{-2}\left(t+1\right)\left(t^{20}\alpha_k^{2}-t^{18}\alpha_k^{2}+t^{17}\alpha_k^{2}-t^{16}\alpha_k^{2}-t^{14}\alpha_k^{2}+2t^{13}\alpha_k^{2}-2t^{12}\alpha_k^{2}+2t^{11}\alpha_k^{2}-3t^{10}\alpha_k^{2}\right. \\ 
 \\ 
 &&\left. \hspace{.1 in}+2t^{9}\alpha_k^{2}-t^{8}\alpha_k^{2}+t^{7}\alpha_k^{2}-t^{6}\alpha_k^{2}+2t^{5}\alpha_k^{2}-t^{4}\alpha_k^{2}+2t^{10}\alpha_k-t^{8}\alpha_k+t^{7}\alpha_k\right. \\ 
 \\ 
 &&\left. \hspace{.1 in}-2t^{5}\alpha_k+t^{3}\alpha_k-t^{2}\alpha_k+1\right)\\ 
\\ 
3/25&5&-\left(t+1\right)\left(t^{32}\alpha_k^{2}-2t^{31}\alpha_k^{2}+t^{30}\alpha_k^{2}-t^{29}\alpha_k^{2}+t^{28}\alpha_k^{2}-2t^{27}\alpha_k^{2}+3t^{26}\alpha_k^{2}-2t^{25}\alpha_k^{2}+2t^{24}\alpha_k^{2}\right. \\ 
 \\ 
 &&\left. \hspace{.1 in}-2t^{23}\alpha_k^{2}+t^{22}\alpha_k^{2}+t^{20}\alpha_k^{2}-t^{19}\alpha_k^{2}+t^{18}\alpha_k^{2}-t^{16}\alpha_k^{2}+2t^{22}\alpha_k-4t^{21}\alpha_k+2t^{20}\alpha_k\right. \\ 
 \\ 
 &&\left. \hspace{.1 in}-2t^{19}\alpha_k+2t^{18}\alpha_k-4t^{17}\alpha_k+6t^{16}\alpha_k-4t^{15}\alpha_k+3t^{14}\alpha_k-3t^{13}\alpha_k+2t^{12}\alpha_k-2t^{11}\alpha_k\right. \\ 
 \\ 
 &&\left. \hspace{.1 in}+2t^{10}\alpha_k-t^{9}\alpha_k+t^{8}\alpha_k+t^{12}-2t^{11}+t^{10}-t^{9}+t^{8}-2t^{7}+3t^{6}-2t^{5}+t^{4}-t^{3}+t^{2}-2t\right. \\ 
 \\ 
 &&\left. \hspace{.1 in}+1\right)\\ 
\end{array}$$
 \end{center}
\renewcommand{\arraystretch}{.5}
 \begin{center}$$\begin{array}{ccl}
p/q &\ell&f(t)\\ 
 \\ 
\hline 
\\ 
7/15&5&\left(t+1\right)\left(t^{20}\alpha_k^{2}-t^{18}\alpha_k^{2}+t^{17}\alpha_k^{2}-t^{16}\alpha_k^{2}-t^{14}\alpha_k^{2}+2t^{13}\alpha_k^{2}-2t^{12}\alpha_k^{2}+2t^{11}\alpha_k^{2}-3t^{10}\alpha_k^{2}\right. \\ 
 \\ 
 &&\left. \hspace{.1 in}+2t^{9}\alpha_k^{2}-t^{8}\alpha_k^{2}+t^{7}\alpha_k^{2}-t^{6}\alpha_k^{2}+2t^{5}\alpha_k^{2}-t^{4}\alpha_k^{2}+2t^{10}\alpha_k-t^{8}\alpha_k+t^{7}\alpha_k\right. \\ 
 \\ 
 &&\left. \hspace{.1 in}-2t^{5}\alpha_k+t^{3}\alpha_k-t^{2}\alpha_k+1\right)\\ 
\\ 
7/25&5&t^{-6}\left(t+1\right)\left(t^{28}\alpha_k^{2}-2t^{26}\alpha_k^{2}+2t^{24}\alpha_k^{2}-t^{23}\alpha_k^{2}-t^{22}\alpha_k^{2}+3t^{21}\alpha_k^{2}-t^{20}\alpha_k^{2}-t^{19}\alpha_k^{2}\right. \\ 
 \\ 
 &&\left. \hspace{.1 in}-2t^{18}\alpha_k^{2}+2t^{17}\alpha_k^{2}+2t^{15}\alpha_k^{2}-4t^{14}\alpha_k^{2}+3t^{13}\alpha_k^{2}-t^{12}\alpha_k^{2}-t^{11}\alpha_k^{2}-t^{10}\alpha_k^{2}+3t^{9}\alpha_k^{2}\right. \\ 
 \\ 
 &&\left. \hspace{.1 in}-t^{8}\alpha_k^{2}+2t^{18}\alpha_k-4t^{16}\alpha_k+5t^{14}\alpha_k-t^{13}\alpha_k-3t^{12}\alpha_k+3t^{11}\alpha_k+t^{10}\alpha_k-4t^{9}\alpha_k\right. \\ 
 \\ 
 &&\left. \hspace{.1 in}-t^{8}\alpha_k+3t^{7}\alpha_k+t^{6}\alpha_k-t^{5}\alpha_k-t^{4}\alpha_k+t^{8}-2t^{6}+3t^{4}-2t^{2}+1\right)\\ 
\\ 
7/45&5&-\left(t+1\right)\left(t^{32}\alpha_k^{2}-3t^{31}\alpha_k^{2}+t^{30}\alpha_k^{2}+t^{29}\alpha_k^{2}+t^{28}\alpha_k^{2}-3t^{27}\alpha_k^{2}+4t^{26}\alpha_k^{2}-2t^{25}\alpha_k^{2}\right. \\ 
 \\ 
 &&\left. \hspace{.1 in}-2t^{23}\alpha_k^{2}+2t^{22}\alpha_k^{2}+t^{21}\alpha_k^{2}+t^{20}\alpha_k^{2}-3t^{19}\alpha_k^{2}+t^{18}\alpha_k^{2}+t^{17}\alpha_k^{2}-2t^{16}\alpha_k^{2}+2t^{14}\alpha_k^{2}\right. \\ 
 \\ 
 &&\left. \hspace{.1 in}-t^{12}\alpha_k^{2}+2t^{22}\alpha_k-6t^{21}\alpha_k+2t^{20}\alpha_k+2t^{19}\alpha_k+2t^{18}\alpha_k-6t^{17}\alpha_k+7t^{16}\alpha_k-5t^{15}\alpha_k\right. \\ 
 \\ 
 &&\left. \hspace{.1 in}+t^{14}\alpha_k-t^{13}\alpha_k+3t^{12}\alpha_k-2t^{11}\alpha_k+3t^{10}\alpha_k-3t^{9}\alpha_k-t^{8}\alpha_k+t^{7}\alpha_k+t^{6}\alpha_k\right. \\ 
 \\ 
 &&\left. \hspace{.1 in}+t^{12}-3t^{11}+t^{10}+t^{9}+t^{8}-3t^{7}+3t^{6}-3t^{5}+t^{4}+t^{3}+t^{2}-3t+1\right)\\ 
\\ 
7/55&5&-t^{-2}\left(t+1\right)\left(t^{32}\alpha_k^{2}-2t^{31}\alpha_k^{2}+t^{30}\alpha_k^{2}-t^{29}\alpha_k^{2}+t^{28}\alpha_k^{2}-2t^{27}\alpha_k^{2}+3t^{26}\alpha_k^{2}-2t^{25}\alpha_k^{2}\right. \\ 
 \\ 
 &&\left. \hspace{.1 in}+2t^{24}\alpha_k^{2}-2t^{23}\alpha_k^{2}+t^{22}\alpha_k^{2}+t^{20}\alpha_k^{2}-t^{19}\alpha_k^{2}+t^{18}\alpha_k^{2}-t^{16}\alpha_k^{2}+2t^{22}\alpha_k-4t^{21}\alpha_k\right. \\ 
 \\ 
 &&\left. \hspace{.1 in}+2t^{20}\alpha_k-2t^{19}\alpha_k+2t^{18}\alpha_k-4t^{17}\alpha_k+6t^{16}\alpha_k-4t^{15}\alpha_k+3t^{14}\alpha_k-3t^{13}\alpha_k+2t^{12}\alpha_k\right. \\ 
 \\ 
 &&\left. \hspace{.1 in}-2t^{11}\alpha_k+2t^{10}\alpha_k-t^{9}\alpha_k+t^{8}\alpha_k+t^{12}-2t^{11}+t^{10}-t^{9}+t^{8}-2t^{7}+3t^{6}-2t^{5}+t^{4}\right. \\ 
 \\ 
 &&\left. \hspace{.1 in}-t^{3}+t^{2}-2t+1\right)\\ 
\\ 
7/65&5&-t^{-2}\left(t+1\right)\left(t^{36}\alpha_k^{2}-2t^{33}\alpha_k^{2}+t^{32}\alpha_k^{2}-2t^{31}\alpha_k^{2}+t^{30}\alpha_k^{2}-t^{29}\alpha_k^{2}+3t^{28}\alpha_k^{2}-2t^{27}\alpha_k^{2}\right. \\ 
 \\ 
 &&\left. \hspace{.1 in}+t^{26}\alpha_k^{2}-2t^{25}\alpha_k^{2}+2t^{24}\alpha_k^{2}+t^{22}\alpha_k^{2}+t^{20}\alpha_k^{2}-t^{19}\alpha_k^{2}-t^{18}\alpha_k^{2}+2t^{26}\alpha_k-4t^{23}\alpha_k\right. \\ 
 \\ 
 &&\left. \hspace{.1 in}+2t^{22}\alpha_k-4t^{21}\alpha_k+2t^{20}\alpha_k-2t^{19}\alpha_k+6t^{18}\alpha_k-3t^{17}\alpha_k+2t^{16}\alpha_k-4t^{15}\alpha_k+3t^{14}\alpha_k\right. \\ 
 \\ 
 &&\left. \hspace{.1 in}-2t^{13}\alpha_k+t^{12}\alpha_k+2t^{10}\alpha_k-t^{9}\alpha_k+t^{16}-2t^{13}+t^{12}-2t^{11}+t^{10}-t^{9}+3t^{8}-t^{7}+t^{6}\right. \\ 
 \\ 
 &&\left. \hspace{.1 in}-2t^{5}+t^{4}-2t^{3}+1\right)\\ 
\\ 
{\bf 7/75}&5&t^{-2}\left(t+1\right)\left(t^{40}\alpha_k^{2}+t^{39}\alpha_k^{2}-t^{38}\alpha_k^{2}-t^{36}\alpha_k^{2}-2t^{34}\alpha_k^{2}+2t^{33}\alpha_k^{2}-t^{32}\alpha_k^{2}+2t^{31}\alpha_k^{2}\right. \\ 
 \\ 
 &&\left. \hspace{.1 in}-3t^{30}\alpha_k^{2}+t^{29}\alpha_k^{2}-t^{28}\alpha_k^{2}+2t^{27}\alpha_k^{2}-t^{26}\alpha_k^{2}+2t^{25}\alpha_k^{2}-t^{22}\alpha_k^{2}+2t^{30}\alpha_k+2t^{29}\alpha_k\right. \\ 
 \\ 
 &&\left. \hspace{.1 in}-2t^{28}\alpha_k-2t^{26}\alpha_k-4t^{24}\alpha_k+4t^{23}\alpha_k-2t^{22}\alpha_k+4t^{21}\alpha_k-4t^{20}\alpha_k+3t^{19}\alpha_k-2t^{18}\alpha_k\right. \\ 
 \\ 
 &&\left. \hspace{.1 in}+4t^{17}\alpha_k-3t^{16}\alpha_k+2t^{15}\alpha_k-t^{14}\alpha_k-2t^{12}\alpha_k+t^{11}\alpha_k+t^{20}+t^{19}-t^{18}-t^{16}-2t^{14}\right. \\ 
 \\ 
 &&\left. \hspace{.1 in}+2t^{13}-t^{12}+2t^{11}-t^{10}+2t^{9}-t^{8}+2t^{7}-2t^{6}-t^{4}-t^{2}+t+1\right)\\ 
\end{array}$$
 \end{center}
\renewcommand{\arraystretch}{.5}
 \begin{center}$$\begin{array}{ccl}
p/q &\ell&f(t)\\ 
 \\ 
\hline 
\\ 
3/7&7&\left(t+1\right)\left(t^{42}\alpha_k^{3}-2t^{40}\alpha_k^{3}+3t^{39}\alpha_k^{3}-4t^{38}\alpha_k^{3}+5t^{37}\alpha_k^{3}-6t^{36}\alpha_k^{3}+6t^{35}\alpha_k^{3}-8t^{34}\alpha_k^{3}\right. \\ 
 \\ 
 &&\left. \hspace{.1 in}+11t^{33}\alpha_k^{3}-13t^{32}\alpha_k^{3}+15t^{31}\alpha_k^{3}-17t^{30}\alpha_k^{3}+19t^{29}\alpha_k^{3}-23t^{28}\alpha_k^{3}+23t^{27}\alpha_k^{3}-21t^{26}\alpha_k^{3}\right. \\ 
 \\ 
 &&\left. \hspace{.1 in}+21t^{25}\alpha_k^{3}-21t^{24}\alpha_k^{3}+21t^{23}\alpha_k^{3}-21t^{22}\alpha_k^{3}+23t^{21}\alpha_k^{3}-21t^{20}\alpha_k^{3}+17t^{19}\alpha_k^{3}-15t^{18}\alpha_k^{3}\right. \\ 
 \\ 
 &&\left. \hspace{.1 in}+13t^{17}\alpha_k^{3}-11t^{16}\alpha_k^{3}+9t^{15}\alpha_k^{3}-6t^{14}\alpha_k^{3}+5t^{13}\alpha_k^{3}-5t^{12}\alpha_k^{3}+4t^{11}\alpha_k^{3}-3t^{10}\alpha_k^{3}\right. \\ 
 \\ 
 &&\left. \hspace{.1 in}+2t^{9}\alpha_k^{3}-t^{8}\alpha_k^{3}-t^{7}\alpha_k^{3}+t^{6}\alpha_k^{3}+3t^{28}\alpha_k^{2}-4t^{26}\alpha_k^{2}+6t^{25}\alpha_k^{2}-7t^{24}\alpha_k^{2}\right. \\ 
 \\ 
 &&\left. \hspace{.1 in}+7t^{23}\alpha_k^{2}-7t^{22}\alpha_k^{2}+4t^{21}\alpha_k^{2}-7t^{20}\alpha_k^{2}+11t^{19}\alpha_k^{2}-13t^{18}\alpha_k^{2}+14t^{17}\alpha_k^{2}-14t^{16}\alpha_k^{2}\right. \\ 
 \\ 
 &&\left. \hspace{.1 in}+14t^{15}\alpha_k^{2}-17t^{14}\alpha_k^{2}+14t^{13}\alpha_k^{2}-10t^{12}\alpha_k^{2}+8t^{11}\alpha_k^{2}-7t^{10}\alpha_k^{2}+7t^{9}\alpha_k^{2}-7t^{8}\alpha_k^{2}\right. \\ 
 \\ 
 &&\left. \hspace{.1 in}+10t^{7}\alpha_k^{2}-7t^{6}\alpha_k^{2}+3t^{5}\alpha_k^{2}-t^{4}\alpha_k^{2}+3t^{14}\alpha_k-2t^{12}\alpha_k+3t^{11}\alpha_k-3t^{10}\alpha_k\right. \\ 
 \\ 
 &&\left. \hspace{.1 in}+2t^{9}\alpha_k-3t^{7}\alpha_k+2t^{5}\alpha_k-3t^{4}\alpha_k+3t^{3}\alpha_k-2t^{2}\alpha_k+1\right)\\ 
\\ 
3/35&7&-t^{-3}\left(t+1\right)\left(t^{72}\alpha_k^{3}-t^{71}\alpha_k^{3}-t^{70}\alpha_k^{3}+2t^{69}\alpha_k^{3}-3t^{68}\alpha_k^{3}+4t^{67}\alpha_k^{3}-5t^{66}\alpha_k^{3}\right. \\ 
 \\ 
 &&\left. \hspace{.1 in}+5t^{65}\alpha_k^{3}-6t^{64}\alpha_k^{3}+9t^{63}\alpha_k^{3}-11t^{62}\alpha_k^{3}+13t^{61}\alpha_k^{3}-15t^{60}\alpha_k^{3}+17t^{59}\alpha_k^{3}-21t^{58}\alpha_k^{3}\right. \\ 
 \\ 
 &&\left. \hspace{.1 in}+23t^{57}\alpha_k^{3}-21t^{56}\alpha_k^{3}+21t^{55}\alpha_k^{3}-21t^{54}\alpha_k^{3}+21t^{53}\alpha_k^{3}-21t^{52}\alpha_k^{3}+23t^{51}\alpha_k^{3}-23t^{50}\alpha_k^{3}\right. \\ 
 \\ 
 &&\left. \hspace{.1 in}+19t^{49}\alpha_k^{3}-17t^{48}\alpha_k^{3}+15t^{47}\alpha_k^{3}-13t^{46}\alpha_k^{3}+11t^{45}\alpha_k^{3}-8t^{44}\alpha_k^{3}+6t^{43}\alpha_k^{3}-6t^{42}\alpha_k^{3}\right. \\ 
 \\ 
 &&\left. \hspace{.1 in}+5t^{41}\alpha_k^{3}-4t^{40}\alpha_k^{3}+3t^{39}\alpha_k^{3}-2t^{38}\alpha_k^{3}+t^{36}\alpha_k^{3}+3t^{58}\alpha_k^{2}-3t^{57}\alpha_k^{2}-3t^{56}\alpha_k^{2}\right. \\ 
 \\ 
 &&\left. \hspace{.1 in}+6t^{55}\alpha_k^{2}-9t^{54}\alpha_k^{2}+12t^{53}\alpha_k^{2}-15t^{52}\alpha_k^{2}+15t^{51}\alpha_k^{2}-18t^{50}\alpha_k^{2}+27t^{49}\alpha_k^{2}-33t^{48}\alpha_k^{2}\right. \\ 
 \\ 
 &&\left. \hspace{.1 in}+39t^{47}\alpha_k^{2}-44t^{46}\alpha_k^{2}+48t^{45}\alpha_k^{2}-56t^{44}\alpha_k^{2}+59t^{43}\alpha_k^{2}-56t^{42}\alpha_k^{2}+56t^{41}\alpha_k^{2}-56t^{40}\alpha_k^{2}\right. \\ 
 \\ 
 &&\left. \hspace{.1 in}+55t^{39}\alpha_k^{2}-53t^{38}\alpha_k^{2}+55t^{37}\alpha_k^{2}-52t^{36}\alpha_k^{2}+43t^{35}\alpha_k^{2}-37t^{34}\alpha_k^{2}+31t^{33}\alpha_k^{2}-26t^{32}\alpha_k^{2}\right. \\ 
 \\ 
 &&\left. \hspace{.1 in}+22t^{31}\alpha_k^{2}-17t^{30}\alpha_k^{2}+14t^{29}\alpha_k^{2}-11t^{28}\alpha_k^{2}+8t^{27}\alpha_k^{2}-5t^{26}\alpha_k^{2}+3t^{25}\alpha_k^{2}-2t^{24}\alpha_k^{2}\right. \\ 
 \\ 
 &&\left. \hspace{.1 in}+3t^{44}\alpha_k-3t^{43}\alpha_k-3t^{42}\alpha_k+6t^{41}\alpha_k-9t^{40}\alpha_k+12t^{39}\alpha_k-15t^{38}\alpha_k+15t^{37}\alpha_k\right. \\ 
 \\ 
 &&\left. \hspace{.1 in}-18t^{36}\alpha_k+27t^{35}\alpha_k-33t^{34}\alpha_k+39t^{33}\alpha_k-43t^{32}\alpha_k+45t^{31}\alpha_k-49t^{30}\alpha_k+49t^{29}\alpha_k\right. \\ 
 \\ 
 &&\left. \hspace{.1 in}-49t^{28}\alpha_k+49t^{27}\alpha_k-49t^{26}\alpha_k+47t^{25}\alpha_k-43t^{24}\alpha_k+41t^{23}\alpha_k-35t^{22}\alpha_k+29t^{21}\alpha_k\right. \\ 
 \\ 
 &&\left. \hspace{.1 in}-25t^{20}\alpha_k+20t^{19}\alpha_k-16t^{18}\alpha_k+13t^{17}\alpha_k-10t^{16}\alpha_k+7t^{15}\alpha_k-4t^{14}\alpha_k+3t^{13}\alpha_k\right. \\ 
 \\ 
 &&\left. \hspace{.1 in}-t^{12}\alpha_k+t^{30}-t^{29}-t^{28}+2t^{27}-3t^{26}+4t^{25}-5t^{24}+5t^{23}-6t^{22}+9t^{21}-11t^{20}+13t^{19}\right. \\ 
 \\ 
 &&\left. \hspace{.1 in}-14t^{18}+14t^{17}-14t^{16}+13t^{15}-14t^{14}+14t^{13}-14t^{12}+13t^{11}-11t^{10}+9t^{9}-6t^{8}+5t^{7}\right. \\ 
 \\ 
 &&\left. \hspace{.1 in}-5t^{6}+4t^{5}-3t^{4}+2t^{3}-t^{2}-t+1\right)\\ 
\end{array}$$
 \end{center}
\renewcommand{\arraystretch}{.43}
 \begin{center}$$\begin{array}{ccl}
p/q &\ell&f(t)\\ 
 \\ 
\hline 
\\ 
5/7&7&\left(t+1\right)\left(t^{42}\alpha_k^{3}+2t^{41}\alpha_k^{3}-4t^{40}\alpha_k^{3}+2t^{39}\alpha_k^{3}-3t^{38}\alpha_k^{3}+5t^{37}\alpha_k^{3}-6t^{36}\alpha_k^{3}+6t^{35}\alpha_k^{3}\right. \\ 
 \\ 
 &&\left. \hspace{.1 in}-10t^{34}\alpha_k^{3}+13t^{33}\alpha_k^{3}-12t^{32}\alpha_k^{3}+14t^{31}\alpha_k^{3}-17t^{30}\alpha_k^{3}+19t^{29}\alpha_k^{3}-23t^{28}\alpha_k^{3}+19t^{27}\alpha_k^{3}\right. \\ 
 \\ 
 &&\left. \hspace{.1 in}-17t^{26}\alpha_k^{3}+23t^{25}\alpha_k^{3}-23t^{24}\alpha_k^{3}+21t^{23}\alpha_k^{3}-21t^{22}\alpha_k^{3}+23t^{21}\alpha_k^{3}-17t^{20}\alpha_k^{3}+13t^{19}\alpha_k^{3}\right. \\ 
 \\ 
 &&\left. \hspace{.1 in}-17t^{18}\alpha_k^{3}+15t^{17}\alpha_k^{3}-11t^{16}\alpha_k^{3}+9t^{15}\alpha_k^{3}-6t^{14}\alpha_k^{3}+7t^{13}\alpha_k^{3}-7t^{12}\alpha_k^{3}+3t^{11}\alpha_k^{3}\right. \\ 
 \\ 
 &&\left. \hspace{.1 in}-2t^{10}\alpha_k^{3}+2t^{9}\alpha_k^{3}-t^{8}\alpha_k^{3}-t^{7}\alpha_k^{3}-t^{6}\alpha_k^{3}+2t^{5}\alpha_k^{3}+t^{4}\alpha_k^{3}-t^{3}\alpha_k^{3}\right. \\ 
 \\ 
 &&\left. \hspace{.1 in}+3t^{28}\alpha_k^{2}+4t^{27}\alpha_k^{2}-7t^{26}\alpha_k^{2}+3t^{25}\alpha_k^{2}-4t^{24}\alpha_k^{2}+6t^{23}\alpha_k^{2}-7t^{22}\alpha_k^{2}+4t^{21}\alpha_k^{2}\right. \\ 
 \\ 
 &&\left. \hspace{.1 in}-11t^{20}\alpha_k^{2}+14t^{19}\alpha_k^{2}-10t^{18}\alpha_k^{2}+11t^{17}\alpha_k^{2}-13t^{16}\alpha_k^{2}+14t^{15}\alpha_k^{2}-17t^{14}\alpha_k^{2}+10t^{13}\alpha_k^{2}\right. \\ 
 \\ 
 &&\left. \hspace{.1 in}-7t^{12}\alpha_k^{2}+11t^{11}\alpha_k^{2}-10t^{10}\alpha_k^{2}+8t^{9}\alpha_k^{2}-7t^{8}\alpha_k^{2}+10t^{7}\alpha_k^{2}-3t^{6}\alpha_k^{2}-4t^{4}\alpha_k^{2}\right. \\ 
 \\ 
 &&\left. \hspace{.1 in}+3t^{3}\alpha_k^{2}-t^{2}\alpha_k^{2}+3t^{14}\alpha_k+2t^{13}\alpha_k-3t^{12}\alpha_k+3t^{9}\alpha_k-2t^{8}\alpha_k-3t^{7}\alpha_k\right. \\ 
 \\ 
 &&\left. \hspace{.1 in}-2t^{6}\alpha_k+3t^{5}\alpha_k-3t^{2}\alpha_k+2t\alpha_k+1\right)\\ 
\\ 
5/21&7&-\left(t+1\right)\left(t^{48}\alpha_k^{3}-4t^{47}\alpha_k^{3}+6t^{46}\alpha_k^{3}-7t^{45}\alpha_k^{3}+8t^{44}\alpha_k^{3}-9t^{43}\alpha_k^{3}+10t^{42}\alpha_k^{3}\right. \\ 
 \\ 
 &&\left. \hspace{.1 in}-12t^{41}\alpha_k^{3}+16t^{40}\alpha_k^{3}-19t^{39}\alpha_k^{3}+21t^{38}\alpha_k^{3}-23t^{37}\alpha_k^{3}+25t^{36}\alpha_k^{3}-27t^{35}\alpha_k^{3}+27t^{34}\alpha_k^{3}\right. \\ 
 \\ 
 &&\left. \hspace{.1 in}-23t^{33}\alpha_k^{3}+21t^{32}\alpha_k^{3}-21t^{31}\alpha_k^{3}+21t^{30}\alpha_k^{3}-21t^{29}\alpha_k^{3}+21t^{28}\alpha_k^{3}-19t^{27}\alpha_k^{3}+13t^{26}\alpha_k^{3}\right. \\ 
 \\ 
 &&\left. \hspace{.1 in}-9t^{25}\alpha_k^{3}+7t^{24}\alpha_k^{3}-5t^{23}\alpha_k^{3}+3t^{22}\alpha_k^{3}-t^{21}\alpha_k^{3}-t^{19}\alpha_k^{3}+t^{18}\alpha_k^{3}-t^{16}\alpha_k^{3}\right. \\ 
 \\ 
 &&\left. \hspace{.1 in}+2t^{15}\alpha_k^{3}-3t^{14}\alpha_k^{3}+3t^{13}\alpha_k^{3}-t^{12}\alpha_k^{3}+3t^{34}\alpha_k^{2}-12t^{33}\alpha_k^{2}+18t^{32}\alpha_k^{2}-21t^{31}\alpha_k^{2}\right. \\ 
 \\ 
 &&\left. \hspace{.1 in}+22t^{30}\alpha_k^{2}-22t^{29}\alpha_k^{2}+21t^{28}\alpha_k^{2}-24t^{27}\alpha_k^{2}+33t^{26}\alpha_k^{2}-39t^{25}\alpha_k^{2}+42t^{24}\alpha_k^{2}-43t^{23}\alpha_k^{2}\right. \\ 
 \\ 
 &&\left. \hspace{.1 in}+43t^{22}\alpha_k^{2}-42t^{21}\alpha_k^{2}+39t^{20}\alpha_k^{2}-30t^{19}\alpha_k^{2}+24t^{18}\alpha_k^{2}-21t^{17}\alpha_k^{2}+20t^{16}\alpha_k^{2}-20t^{15}\alpha_k^{2}\right. \\ 
 \\ 
 &&\left. \hspace{.1 in}+21t^{14}\alpha_k^{2}-18t^{13}\alpha_k^{2}+9t^{12}\alpha_k^{2}-3t^{11}\alpha_k^{2}+t^{9}\alpha_k^{2}-t^{8}\alpha_k^{2}+3t^{20}\alpha_k-12t^{19}\alpha_k\right. \\ 
 \\ 
 &&\left. \hspace{.1 in}+18t^{18}\alpha_k-21t^{17}\alpha_k+20t^{16}\alpha_k-17t^{15}\alpha_k+12t^{14}\alpha_k-12t^{13}\alpha_k+17t^{12}\alpha_k-20t^{11}\alpha_k\right. \\ 
 \\ 
 &&\left. \hspace{.1 in}+21t^{10}\alpha_k-20t^{9}\alpha_k+17t^{8}\alpha_k-12t^{7}\alpha_k+9t^{6}\alpha_k-5t^{5}\alpha_k+2t^{4}\alpha_k+t^{6}\right. \\ 
 \\ 
 &&\left. \hspace{.1 in}-4t^{5}+6t^{4}-7t^{3}+6t^{2}-4t+1\right)\\ 
\\ 
5/49&7&-t^{-3}\left(t+1\right)\left(t^{66}\alpha_k^{3}-3t^{65}\alpha_k^{3}+3t^{64}\alpha_k^{3}-2t^{63}\alpha_k^{3}+t^{62}\alpha_k^{3}-t^{60}\alpha_k^{3}+t^{59}\alpha_k^{3}\right. \\ 
 \\ 
 &&\left. \hspace{.1 in}+t^{57}\alpha_k^{3}-3t^{56}\alpha_k^{3}+5t^{55}\alpha_k^{3}-7t^{54}\alpha_k^{3}+9t^{53}\alpha_k^{3}-13t^{52}\alpha_k^{3}+19t^{51}\alpha_k^{3}-21t^{50}\alpha_k^{3}\right. \\ 
 \\ 
 &&\left. \hspace{.1 in}+21t^{49}\alpha_k^{3}-21t^{48}\alpha_k^{3}+21t^{47}\alpha_k^{3}-21t^{46}\alpha_k^{3}+23t^{45}\alpha_k^{3}-27t^{44}\alpha_k^{3}+27t^{43}\alpha_k^{3}-25t^{42}\alpha_k^{3}\right. \\ 
 \\ 
 &&\left. \hspace{.1 in}+23t^{41}\alpha_k^{3}-21t^{40}\alpha_k^{3}+19t^{39}\alpha_k^{3}-16t^{38}\alpha_k^{3}+12t^{37}\alpha_k^{3}-10t^{36}\alpha_k^{3}+9t^{35}\alpha_k^{3}-8t^{34}\alpha_k^{3}\right. \\ 
 \\ 
 &&\left. \hspace{.1 in}+7t^{33}\alpha_k^{3}-6t^{32}\alpha_k^{3}+4t^{31}\alpha_k^{3}-t^{30}\alpha_k^{3}+3t^{52}\alpha_k^{2}-9t^{51}\alpha_k^{2}+9t^{50}\alpha_k^{2}-6t^{49}\alpha_k^{2}\right. \\ 
 \\ 
 &&\left. \hspace{.1 in}+3t^{48}\alpha_k^{2}-3t^{46}\alpha_k^{2}+3t^{45}\alpha_k^{2}+3t^{43}\alpha_k^{2}-10t^{42}\alpha_k^{2}+16t^{41}\alpha_k^{2}-21t^{40}\alpha_k^{2}+24t^{39}\alpha_k^{2}\right. \\ 
 \\ 
 &&\left. \hspace{.1 in}-30t^{38}\alpha_k^{2}+39t^{37}\alpha_k^{2}-42t^{36}\alpha_k^{2}+43t^{35}\alpha_k^{2}-43t^{34}\alpha_k^{2}+42t^{33}\alpha_k^{2}-39t^{32}\alpha_k^{2}+39t^{31}\alpha_k^{2}\right. \\ 
 \\ 
 &&\left. \hspace{.1 in}-42t^{30}\alpha_k^{2}+39t^{29}\alpha_k^{2}-32t^{28}\alpha_k^{2}+26t^{27}\alpha_k^{2}-21t^{26}\alpha_k^{2}+18t^{25}\alpha_k^{2}-15t^{24}\alpha_k^{2}+12t^{23}\alpha_k^{2}\right. \\ 
 \\ 
 &&\left. \hspace{.1 in}-9t^{22}\alpha_k^{2}+5t^{21}\alpha_k^{2}-2t^{20}\alpha_k^{2}+3t^{38}\alpha_k-9t^{37}\alpha_k+9t^{36}\alpha_k-6t^{35}\alpha_k+3t^{34}\alpha_k\right. \\ 
 \\ 
 &&\left. \hspace{.1 in}-3t^{32}\alpha_k+3t^{31}\alpha_k+3t^{29}\alpha_k-11t^{28}\alpha_k+17t^{27}\alpha_k-21t^{26}\alpha_k+21t^{25}\alpha_k-21t^{24}\alpha_k\right. \\ 
 \\ 
 &&\left. \hspace{.1 in}+21t^{23}\alpha_k-21t^{22}\alpha_k+23t^{21}\alpha_k-23t^{20}\alpha_k+21t^{19}\alpha_k-17t^{18}\alpha_k+14t^{17}\alpha_k-12t^{16}\alpha_k\right. \\ 
 \\ 
 &&\left. \hspace{.1 in}+9t^{15}\alpha_k-6t^{14}\alpha_k+3t^{13}\alpha_k-t^{11}\alpha_k+t^{10}\alpha_k+t^{24}-3t^{23}+3t^{22}-2t^{21}+t^{20}\right. \\ 
 \\ 
 &&\left. \hspace{.1 in}-t^{18}+t^{17}+t^{15}-4t^{14}+6t^{13}-7t^{12}+6t^{11}-4t^{10}+t^{9}+t^{7}-t^{6}+t^{4}-2t^{3}+3t^{2}-3t\right. \\ 
 \\ 
 &&\left. \hspace{.1 in}+1\right)\\ 
\end{array}$$
 \end{center}

\renewcommand{\arraystretch}{.5}
 \begin{center}$$\begin{array}{ccl}
p/q &\ell&f(t)\\ 
 \\ 
\hline 
\\
5/63&7&t^{-6}\left(t+1\right)\left(t^{78}\alpha_k^{3}-t^{77}\alpha_k^{3}-2t^{76}\alpha_k^{3}+t^{75}\alpha_k^{3}+t^{74}\alpha_k^{3}+t^{73}\alpha_k^{3}-2t^{72}\alpha_k^{3}\right. \\ 
 \\ 
 &&\left. \hspace{.1 in}+2t^{71}\alpha_k^{3}-3t^{70}\alpha_k^{3}+7t^{69}\alpha_k^{3}-7t^{68}\alpha_k^{3}+6t^{67}\alpha_k^{3}-9t^{66}\alpha_k^{3}+11t^{65}\alpha_k^{3}-15t^{64}\alpha_k^{3}\right. \\ 
 \\ 
 &&\left. \hspace{.1 in}+17t^{63}\alpha_k^{3}-13t^{62}\alpha_k^{3}+17t^{61}\alpha_k^{3}-23t^{60}\alpha_k^{3}+21t^{59}\alpha_k^{3}-21t^{58}\alpha_k^{3}+23t^{57}\alpha_k^{3}-23t^{56}\alpha_k^{3}\right. \\ 
 \\ 
 &&\left. \hspace{.1 in}+17t^{55}\alpha_k^{3}-19t^{54}\alpha_k^{3}+23t^{53}\alpha_k^{3}-19t^{52}\alpha_k^{3}+17t^{51}\alpha_k^{3}-14t^{50}\alpha_k^{3}+12t^{49}\alpha_k^{3}-13t^{48}\alpha_k^{3}\right. \\ 
 \\ 
 &&\left. \hspace{.1 in}+10t^{47}\alpha_k^{3}-6t^{46}\alpha_k^{3}+6t^{45}\alpha_k^{3}-5t^{44}\alpha_k^{3}+3t^{43}\alpha_k^{3}-2t^{42}\alpha_k^{3}+4t^{41}\alpha_k^{3}-2t^{40}\alpha_k^{3}\right. \\ 
 \\ 
 &&\left. \hspace{.1 in}-t^{39}\alpha_k^{3}+3t^{64}\alpha_k^{2}-3t^{63}\alpha_k^{2}-6t^{62}\alpha_k^{2}+3t^{61}\alpha_k^{2}+3t^{60}\alpha_k^{2}+3t^{59}\alpha_k^{2}-6t^{58}\alpha_k^{2}\right. \\ 
 \\ 
 &&\left. \hspace{.1 in}+6t^{57}\alpha_k^{2}-9t^{56}\alpha_k^{2}+21t^{55}\alpha_k^{2}-21t^{54}\alpha_k^{2}+18t^{53}\alpha_k^{2}-27t^{52}\alpha_k^{2}+32t^{51}\alpha_k^{2}-42t^{50}\alpha_k^{2}\right. \\ 
 \\ 
 &&\left. \hspace{.1 in}+47t^{49}\alpha_k^{2}-39t^{48}\alpha_k^{2}+48t^{47}\alpha_k^{2}-59t^{46}\alpha_k^{2}+56t^{45}\alpha_k^{2}-55t^{44}\alpha_k^{2}+59t^{43}\alpha_k^{2}-58t^{42}\alpha_k^{2}\right. \\ 
 \\ 
 &&\left. \hspace{.1 in}+44t^{41}\alpha_k^{2}-47t^{40}\alpha_k^{2}+52t^{39}\alpha_k^{2}-43t^{38}\alpha_k^{2}+38t^{37}\alpha_k^{2}-31t^{36}\alpha_k^{2}+26t^{35}\alpha_k^{2}-25t^{34}\alpha_k^{2}\right. \\ 
 \\ 
 &&\left. \hspace{.1 in}+19t^{33}\alpha_k^{2}-14t^{32}\alpha_k^{2}+11t^{31}\alpha_k^{2}-9t^{30}\alpha_k^{2}+5t^{29}\alpha_k^{2}-3t^{28}\alpha_k^{2}+5t^{27}\alpha_k^{2}-2t^{26}\alpha_k^{2}\right. \\ 
 \\ 
 &&\left. \hspace{.1 in}+3t^{50}\alpha_k-3t^{49}\alpha_k-6t^{48}\alpha_k+3t^{47}\alpha_k+3t^{46}\alpha_k+3t^{45}\alpha_k-6t^{44}\alpha_k+6t^{43}\alpha_k\right. \\ 
 \\ 
 &&\left. \hspace{.1 in}-9t^{42}\alpha_k+21t^{41}\alpha_k-21t^{40}\alpha_k+18t^{39}\alpha_k-27t^{38}\alpha_k+31t^{37}\alpha_k-39t^{36}\alpha_k+43t^{35}\alpha_k\right. \\ 
 \\ 
 &&\left. \hspace{.1 in}-39t^{34}\alpha_k+45t^{33}\alpha_k-49t^{32}\alpha_k+49t^{31}\alpha_k-47t^{30}\alpha_k+49t^{29}\alpha_k-47t^{28}\alpha_k+37t^{27}\alpha_k\right. \\ 
 \\ 
 &&\left. \hspace{.1 in}-37t^{26}\alpha_k+35t^{25}\alpha_k-31t^{24}\alpha_k+28t^{23}\alpha_k-20t^{22}\alpha_k+16t^{21}\alpha_k-14t^{20}\alpha_k+10t^{19}\alpha_k\right. \\ 
 \\ 
 &&\left. \hspace{.1 in}-7t^{18}\alpha_k+6t^{17}\alpha_k-6t^{16}\alpha_k+t^{15}\alpha_k+t^{13}\alpha_k+t^{36}-t^{35}-2t^{34}+t^{33}+t^{32}\right. \\ 
 \\ 
 &&\left. \hspace{.1 in}+t^{31}-2t^{30}+2t^{29}-3t^{28}+7t^{27}-7t^{26}+6t^{25}-9t^{24}+10t^{23}-12t^{22}+13t^{21}-13t^{20}\right. \\ 
 \\ 
 &&\left. \hspace{.1 in}+14t^{19}-13t^{18}+14t^{17}-13t^{16}+13t^{15}-12t^{14}+10t^{13}-9t^{12}+6t^{11}-7t^{10}+7t^{9}\right. \\ 
 \\ 
 &&\left. \hspace{.1 in}-3t^{8}+2t^{7}-2t^{6}+t^{5}+t^{4}+t^{3}-2t^{2}-t+1\right)\\ 
\end{array}$$
 \end{center}

 \end{document}